\input amstex
\documentstyle{amsppt}
\magnification=\magstep1

\pageheight{9.0truein}
\pagewidth{6.5truein}

\NoBlackBoxes

\long\def\ignore#1{#1}

\ignore{
\input xy
\xyoption{matrix}\xyoption{arrow}\xyoption{curve}\xyoption{frame}
\def\edge{\ar@{-}}

\def\threepool{\save[0,0]+(-3,3);[0,2]+(3,-3)**\frm<10pt>{.}\restore}
\def\toplabel#1{\save+<0ex,3ex>\drop{#1}\restore} }

\def\la{\Lambda}
\def\lamod{\Lambda{}\operatorname{-mod}}
\def\Lamod{\Lambda{}\operatorname{-Mod}}
\def\pinf{{\Cal P}^\infty}
\def\pinflamod{{\Cal P}^\infty(\Lambda{}\operatorname{-mod})}

\def\sinflamod{{\Cal S}^\infty(\Lambda{}\operatorname{-mod})}
\def\pinfapprox{$\pinflamod$-ap\-prox\-i\-ma\-tion}
\def\pdim{\operatorname{p\,dim}}
\def\Hom{\operatorname{Hom}}
\def\lfindim{\operatorname{l\,fin\,dim}}
\def\lcycfindim{\operatorname{l\,cyc\,fin\,dim}}
\def\lFindim{\operatorname{l\,Fin\,dim}}

\def\soc{\operatorname{soc}}
\def\NN{\Bbb N}
\def\ZZ{\Bbb Z}
\def\kgami{K\Gamma/I}
\def\A{{\Cal A}}
\def\C{{\Cal C}}
\def\D{{\Cal D}}
\def\filt{\operatorname{filt}}
\def\Filt{\overarrow{\operatorname{filt}}}
\def\that{\widehat{t}}
\def\That{\widehat{T}}
\def\ptil{\widetilde{p}}
\def\qtil{\widetilde{q}}

\def\Alp{{\bf 1}}
\def\AMMH{{\bf 2}}
\def\AuRe{{\bf 3}}
\def\AuSm{{\bf 4}}
\def\AuSmtwo{{\bf 5}}
\def\Bass{{\bf 6}}
\def\Bel{{\bf 7}}
\def\BenGna{{\bf 8}}
\def\Bleone{{\bf 9}}
\def\Bletwo{{\bf 10}}
\def\Bon{{\bf 11}}
\def\BuHZ{{\bf 12}}
\def\BuRi{{\bf 13}}
\def\ChriS{{\bf 14}}
\def\CB{{\bf 15}}
\def\CBser{{\bf 16}}
\def\DoFr{{\bf 17}}
\def\EckSch{{\bf 18}}
\def\Eno{{\bf 19}}
\def\EMDV{{\bf 20}}
\def\ErSk{{\bf 21}}
\def\Ful{{\bf 22}}
\def\Gabriel{{\bf 23}}
\def\Gab{{\bf 24}}
\def\Geitame{{\bf 25}}
\def\Geiser{{\bf 26}}
\def\Geiclan{{\bf 27}}
\def\GePo{{\bf 28}}
\def\HaHZ{{\bf 29}}
\def\pre{{\bf 30}}
\def\dom{{\bf 31}}
\def\convenient{{\bf 32}}
\def\dep{{\bf 33}}
\def\HZSmdnsn{{\bf 34}}
\def\HZSm{{\bf 35}}
\def\HZSmbsrl{{\bf 36}}
\def\IgSmTo{{\bf 37}}
\def\Levy{{\bf 38}}
\def\Mar{{\bf 39}}
\def\McG{{\bf 40}}
\def\RaSa{{\bf 41}}
\def\RSDV{{\bf 42}}
\def\Rin{{\bf 43}}
\def\Rintwo{{\bf 44}}
\def\Ringen{{\bf 45}} 
\def\Rog{{\bf 46}}
\def\SkWa{{\bf 47}}
\def\Sma{{\bf 48}}
\def\Sze{{\bf 49}}
\def\VFCB{{\bf 50}}
\def\WaWa{{\bf 51}}
\def\Xu{{\bf 52}}
\def\Xutwo{{\bf 53}}

\topmatter

\title The phantom menace in representation theory
\endtitle

\author Birge Huisgen-Zimmermann \endauthor

\address Department of Mathematics, University of California, Santa Barbara,
CA 93106, USA\endaddress
\email birge\@math.ucsb.edu\endemail

\thanks While writing this paper, the author was partially supported by a
grant from the National Science Foundation \endthanks

\endtopmatter

\document

\head 1. Introduction and prerequisites\endhead

Our principal goal in this overview is to explain and motivate the concept of
a {\it phantom\/} in the representation theory of finite dimensional
algebras.  In particular, we will exhibit the key role of phantoms towards
understanding how a full subcategory $\A$ of 
$\lamod$ is embedded into $\lamod$, in terms of maps leaving or entering
$\A$. (Here $\la$ is a finite dimensional algebra over a field $K$, and
$\lamod$ denotes the category of all finitely generated left $\la$-modules.) 
Very roughly speaking, phantoms serve the following dual purpose in this
connection: On one hand, they represent an effective tool for systematically
tackling the question of whether or not
$\A$ is functorially finite in $\lamod$, or  --  in less technical terms  -- 
whether all objects in $\lamod$ have best (right/left) approximations in
$\A$, in a sense to be made precise momentarily.  On the other hand, the {\it
effective\/}
$\A$-phan\-toms of a given object $X\in\lamod$ capture -- within a {\it
minimal\/} frame -- a condensed picture of the relations of those objects in
$\A$ which have nontrivial homomorphisms to
$X$.

For a preview of the central definitions, let $\A$ be closed under finite
direct sums and $\C$ a subcategory of
$\A$.  Recall that
$f \in\Hom_\la(A,X)$ is called a $\C$-{\it ap\-prox\-i\-ma\-tion of\/} $X$
{\it inside\/} $\A$ in case
$A$ belongs to $\A$ and all maps in $\Hom_\la(C,X)$ with $C\in\C$ factor
through $f$ \cite{\HaHZ}. In case $\C=\A$, we re-encounter the classical
right $\A$-ap\-prox\-i\-ma\-tions of $X$ as introduced by Auslander-Smal\o\
\cite{\AuSm} and, independently, Enochs \cite{\Eno}.  Whenever such right
approximations exist, there is a unique candidate of minimal $K$-dimension,
which thus serves as `best' (right) approximation of $X$ in
$\A$.  (The concept of left $\A$-approxiamation is dual; we will suppress the
qualifier `right' in the sequel, since we will leave dualization of our
notions and results to the reader.)  Note, however, that
$\C$-ap\-prox\-i\-ma\-tions of $X$ inside
$\A$ need not exist in general.  One can, a priori, rely on their
availability only when $\C$ is finite. So if
$\C$ is countable, for instance, say $\C=\{C_n\mid n\in\NN\}$, it is natural
to consider $\{C_1,\dots,C_n\}$-ap\-prox\-i\-ma\-tions of $X$ inside
$\A$, and to explore whether they `converge' to an object in $\A$.  In
general they will not, but rather keep  growing in dimension.  So the
question arises whether one should drop the requirement that the expected
limits be finitely generated.  Of course, if we do not insist on {\it
finitely generated\/} approximating objects, the factorization problem per se
is trivialized; namely  --  to again refer to our countable test situation 
--  forming the direct sum of a full collection of
$\{C_1,\dots,C_n\}$-ap\-prox\-i\-ma\-tions of $X$ inside $\A$, where $n$
traces $\NN$, will always yield the domain of a homomorphism through which
all maps $C_n\rightarrow X$ can be factored.  Yet, such sums will be highly
redundant with respect to this stipulation as a rule.  The factorization
problem thus re-gains interest if we require that the redundancy be trimmed
off.  The minimality condition which Happel and the author imposed
\cite{\HaHZ} to produce condensed images of the `mapping behavior' of
$\A$ relative to $X$, is as follows: An $\A$-{\it phan\-tom\/} of
$X$ is a module
$H\in\Lamod$ such that each finitely generated submodule $H'$ of $H$ has the
following property: there exists a finite subcategory $\A'\subseteq \A$  -- 
which will usually depend on $H'$  --  such that $H'$ occurs as a subfactor
of {\it every\/} $\A'$-ap\-prox\-i\-ma\-tion of $X$ inside $\A$. On one hand,
this definition is quite loose, in that the class of
$\A$-phan\-toms of $X$ is obviously closed under subfactors and direct limits
(of directed systems).  On the other hand, for {\it finite\/} classes $\A'
\subseteq \A$, one expects a large number of $\A'$-ap\-prox\-i\-ma\-tions of
$X$ inside $\A$, whence the requirement that $H'$ should make an appearance
in all of them places major pressure on $H'$.  The slack in the definition is
useful towards testing whether $X$ has a classical
$\A$-ap\-prox\-i\-ma\-tion, for the following reason:  The existence of
classical approximations obviously forces all
$\A$-phantoms of $X$ to be subfactors of the minimal
$\A$-ap\-prox\-i\-ma\-tion of $X$.  In fact, the class of phantoms contains
infinite dimensional objects  --  or, equivalently, objects of unbounded
finite dimensions  --  if and only if
$X$ fails to have a classical $\A$-ap\-prox\-i\-ma\-tion.  So the easier it
is to construct (infinite dimensional) phantoms, the better, as long as they
are only to serve as indicators for the classical existence problem.  On the
other hand, in case this problem has already been resolved in the negative,
small
 phantoms are no longer of interest.  One is led to investigate a class of
upgraded phantoms that are fairly close in spirit to the classical
approximations, except for usually being infinite dimensional: namely, the
`effective' phantoms.  Given a subcategory $\C \subseteq \A$, a $\C$-{\it
effective phantom\/} $H$ of
$X$ is a direct limit of a directed system of objects in
$\A$ which is a phantom and comes equipped with a homomorphism $H \rightarrow
X$, through which all homomorphisms
$C \rightarrow X$ for $C \in \C$ can be factored. For existence of effective
phantoms, see Section 4.

On the side, we mention that, in algebraic topology, there has long been
interest in {\it phantom maps\/}, namely in maps $f: X \rightarrow Y$ of
CW-complexes with the property that, for each $n \in
\NN$, the restriction of $f$ to the $n$-skeleton $X^n$ of $X$ is
null-homotopic (see \cite{\McG} for a survey and \cite{\ChriS} for a modified
notion of a phantom map).  Around the same time as Happel and the author
introduced the concepts outlined above, Benson and Gnacadja
\cite{\BenGna} defined phantom maps in the context of group representations
as follows:  Given a finite group $G$ and a field $K$ of positive
characteristic, a homomorphism $X \rightarrow Y$ of
$KG$-modules is said to be a phantom map if its restriction to every finitely
generated submodule of $M$ factors through a projective module.  These
phantom maps share a decisive feature with the universal maps accompanying
the effective phantoms we mentioned in the preceding paragraph:  Namely, they
are comparatively well understood on restrictions to finitely generated
subobjects of the domain, while globally, they often display qualitatively
new properties.

Let us add a bit more history.  The approximation theory that gave rise to
the work recorded in this overview has its roots in the discovery that
injective envelopes and projective covers are highly useful auxiliary objects
in the structural analysis of more general modules.  In their seminal paper
of 1953, Eckmann and Schopf
\cite{\EckSch} showed that every module can be embedded as an essential
submodule into an injective module, and in 1962, Gabriel \cite{\Gabriel}
extended their findings to the context of Grothendieck categories.  The more
delicate existence question for projective covers, on the other hand, was
resolved by Bass
\cite{\Bass} in 1960, in his famous Theorem P.  Put in slightly different
terms, these results address best approximations of a module $X$ by an object
from the category of all injective, resp. projective, modules.  The passage
to more general subcategories of modules was performed by Auslander-Smal\o\
\cite{\AuSm} and Enochs \cite{\Eno}, almost simultaneously, in the early
eighties.  The former team was primarily interested in categories of finitely
generated modules over Artin algebras and spoke about right and left
approximations, while Enochs focused on more general settings, such as
coherent or noetherian rings, calling the decisive maps  covers and
envelopes, depending on whether they leave or enter the considered module
category.  The artinian line was first used towards sufficient conditions for
the existence of almost split sequences in subcategories by Auslander and
Smal\o\ 
\cite{\AuSmtwo} in 1981; a decade later, this line obtained an important boost
through its applications to homology and tilting theory, first observed by
Auslander-Reiten
\cite{\AuRe}.  Meanwhile, the latter approach was pursued separately by
Asensio Mayor, Belshoff, Enochs, Guil-Asensio, Mart\'\i nez Hern\'andez, Rada,
Saor\'\i n, del Valle, Xu, and others (see
\cite{\AMMH, \Bel,
\EMDV, \Mar, \RaSa, \RSDV, \Xu, \Xutwo}).  Finally, we mention that in
\cite{\Levy}, Levy independently used minimal approximations from a category
of particulary accessible objects, in order to study modules over pullback
rings.        

In the sequel, our favorite category will be $\pinflamod$, the category of
all finitely generated left
$\la$-modules of finite projective dimension; it is contained in the category
$\pinf(\Lamod)$ having as objects {\it all\/} $\la$-modules of finite
projective dimension.  To place the new concepts (of Section 4) into context,
we will discuss and exemplify (in Sections 2,3) the massive impact of
contravariant finiteness of $\pinflamod$ on the representation theory of
$\la$. The most striking results in this direction are Theorems 2 and 3 of
Section 3 which, in case of existence, identify the minimal
$\pinflamod$-approximations of the simple modules as the basic building
blocks of arbitrary objects in
$\pinf(\Lamod)$  --  and as the indicators of little and big finitistic
dimensions; these theorems combine results of Auslander and Reiten
\cite{\AuRe} with more recent work of Smal\o\ and the author \cite{\HZSm}. 
The resulting program of accessing the structure of objects in
$\pinf(\Lamod)$ can be carried out in an ideally explicit format when
$\la$ is left serial
\cite{\BuHZ}. This first model situation will be presented at the end of
Section 3. The second class of algebras for which the theory developed yields
particularly complete and smooth results is the class of string algebras. In
that case, answers to essentially all homological questions one might pose
can be provided in the form of a finite number of `characteristic'
$\pinflamod$-phan\-toms, one for each of the simple left modules. This is the
content of ongoing joint work of the author with Smal\o\ \cite{\HZSmbsrl}; it
is sketched in Section 5. We include a brief history of the theory of string
algebras, since these algebras have established their role as excellent
display cases of the representation-theoretic methods developed during the
past decades.

Since our primary goal is to familiarize the reader with the type of
information stored in phantoms, and since we consider examples as one of the
keys to an intuitive grasp, there will be major emphasis on explicit
illustrations.  In fact, we start by setting up a sequence of `test examples'
to which we keep returning as the discussion proceeds.

This is a revised and updated version of a survey that appeared in the
Proceedings of the $30$th Symposium on Ring Theory and Representation Theory,
which took place in Nagano, Japan, in 1997.  The author would like to thank
the organizer/editor, Y. Iwanaga, for his permission to incorporate portions
of that earlier publication here.  Moreover, the author is indebted to M.
Saor\'\i n for references to the history of the non-artinian thread of the
subject.

\medskip

{\narrower\noindent {\sl There are people indeed [\dots] for whom all the
things that have a fixed value, assessable by others, fortune, success, high
positions, do not count; what they must have is phantoms.} 
\smallskip}

\rightline{Marcel Proust, {\sl Remembrance of Things Past}}

\subhead Content overview\endsubhead

{\bf 2.} Contravariant finiteness and first examples

{\bf 3.} Homological importance of contravariant finiteness and a model
application of the theory

{\bf 4.} Phantoms. Definitions, existence, and basic properties

{\bf 5.} Phantoms over string algebras  --  another model setting

\subhead Prerequisites\endsubhead

Throughout, $\la=\kgami$ will be a path algebra modulo relations with Jacobson
radical $J$, and the vertex set of $\Gamma$ will be identified with a full set
of primitive idempotents of $\la$. By $\Lamod$ we will denote the category of
all left $\la$-modules and by $\lamod$ the full subcategory of finitely
generated modules. Moreover, $\pinf(\Lamod)$ and $\pinflamod$ will be the
subcategories of $\Lamod$ and $\lamod$, respectively, having as objects the
modules of finite projective dimension. The suprema of the projective
dimensions attained on $\pinf(\Lamod)$ and $\pinflamod$ will be labeled
$\lFindim\la$ and $\lfindim\la$, respectively.

Given a path $p$ in $K\Gamma$, we will denote by $\text{start}(p)$ and
$\text{end}(p)$ the starting and end points of $p$, respectively, and our
convention for multiplying paths $p,q\in K\Gamma$ is as follows: ``$qp$''
stands for ``$q$ after $p$''.

Our most important  auxiliary devices will be labeled and layered graphs of
finite dimensional modules. Since our graphing conventions differ to some
extent from those of other authors (in particular, they are akin to, but not
the same as, the module diagrams studied by Alperin \cite{\Alp} and Fuller
\cite{\Ful}), we include an informal description of our graphs for the
convenience of the reader.

The graphs we use are based on sequences of `top elements' of a module $M$
which are $K$-linearly independent modulo $JM$. Here $x\in M$ is a {\it top
element\/} of $M$ if $x\notin JM$ and $x=ex$ for one of the primitive
idempotents $e$ corresponding to the vertices of $\Gamma$; in this situation
we will also say that $x$ is a top element of {\it type\/} $e$.

Let $\la=\kgami$ be the algebra presented in Example C.1 of Section 2 below.
That the indecomposable projective module $\la e_6$ have  graph 

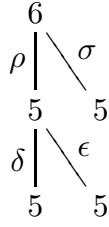
\ignore{
$$\xymatrixcolsep{0.8pc}\xymatrixrowsep{1.5pc}
\xymatrix{ 6\edge[d]_\rho \edge[dr]^\sigma\\ 5\edge[d]_\delta
\edge[dr]^\epsilon &5\\ 5 &5 }$$ }

\noindent with respect to the top element $e_6$ means that $J^3e_6=0$ and
$Je_6/J^2e_6 \cong J^2e_6\cong S_5\oplus S_5$, that $\rho e_6$ and $\sigma
e_6$ are $K$-linearly independent modulo $J^2e_6$, and that $\delta\rho e_6$
and $\epsilon\rho e_6$ are
$K$-linearly independent (modulo $J^3e_6$), while $\delta\sigma e_6=
\epsilon\delta e_6=0$ in
$\la e_6$.

Generally, the entries in each row of a graph of a module $M$ record the
composition factors in the radical layers $J^rM/J^{r+1}M$ in the correct
multiplicities. Whenever we present the graph of an indecomposable projective
module $\la e$, we tacitly assume that the corresponding top element of $\la
e$ is chosen to be $e$. In our first example the choice of top element of
$\la e_6$ does not influence the graph, but it will in other situations. For
instance, the module
$M= (\la e_6\oplus \la e_6)/U$, where $U= \la( \rho e_6, \rho e_6)$, has graph

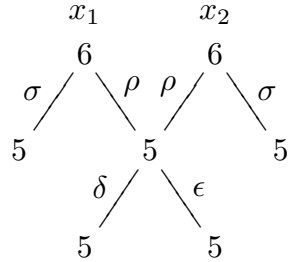
\ignore{
$$\xymatrixcolsep{0.8pc}\xymatrixrowsep{1.5pc}
\xymatrix{
 &6\toplabel{x_1} \edge[dl]_\sigma \edge[dr]^(0.45)\rho &&6\toplabel{x_2}
\edge[dl]_(0.45)\rho \edge[dr]^\sigma\\ 5 &&5\edge[dl]_\delta
\edge[dr]^\epsilon &&5\\
 &5 &&5 }$$ }

\noindent relative to the top elements $x_1=(e_6,0)$ and $x_2= (0,e_6)$, while
its graph relative to the top elements $y_1=x_1+x_2$ and $y_2=x_2$ is

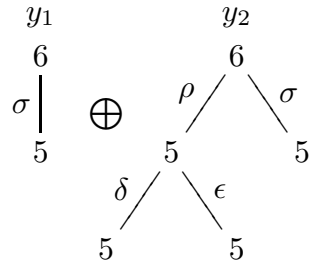
\ignore{
$$\xymatrixcolsep{0.8pc}\xymatrixrowsep{1.5pc}
\xymatrix{ 6\toplabel{y_1} \edge[d]_\sigma &&&6\toplabel{y_2} \edge[dl]_\rho
\edge[dr]^\sigma\\ 5 &\toplabel{\bigoplus} &5\edge[dl]_\delta
\edge[dr]^\epsilon &&5\\
 &5 &&5
 }$$ }

Note that the graph of a module need not determine the module in question, up
to isomorphism. For example, each of the modules $M_k= \la e_6/U_k$, where
$U_k= \la( \sigma- k\rho )e_6$ with $k\in K^*$, has graph

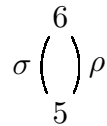
\ignore{
$$\xymatrixcolsep{0.8pc}\xymatrixrowsep{1.5pc}
\xymatrix{ 6\edge@/_/[d]_\sigma \edge@/^/[d]^\rho\\ 5 }$$ }

\noindent with respect to the top element $e_6+ U_k$ of $M_k$, while
$M_{k_1} \not\cong M_{k_2}$ if $k_1\ne k_2$.

For a slightly more involved example, consider the factor module $N= (\la e_6
\oplus \la e_{11}) /V$, where $V$ is the submodule generated by $(\rho,\,
\rho\psi +\delta\rho\psi)$ and $(\epsilon\rho -\sigma,\, 0)$. Relative to the
top elements $x_1= (e_6,0)$ and $x_2= (0,e_{11})$, the module
$N$ has graph

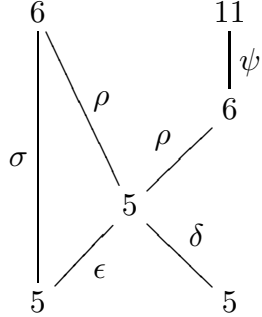
\ignore{
$$\xymatrixcolsep{1.5pc}\xymatrixrowsep{1.5pc}
\xymatrix{ 6 \edge[ddr]^\rho \edge[ddd]_\sigma &&11 \edge[d]^\psi\\
 &&6 \edge[dl]_\rho\\
 &5 \edge[dl]^\epsilon \edge[dr]^\delta\\ 5 &&5 }$$ }

\noindent indicating -- next to the previously discussed information -- that
$\sigma x_1$ is equal to a nonzero scalar multiple of $\epsilon\rho x_1$,
modulo $J^4M=0$, and $\rho x_1$ equals a nonzero scalar multiple of $\rho\psi
x_2$, modulo $J^3M$.

To enlarge the family of objects in $\lamod$ which possess labeled and layered
graphs relative to suitable sequences of top elements, we also allow for
graphs with `{\it pools\/}' along the following model. That a module
$A\in\lamod$ has graph

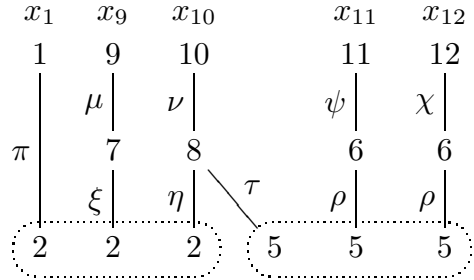
\ignore{
$$\xymatrixcolsep{1.0pc}\xymatrixrowsep{1.5pc}
\xymatrix{ 1\toplabel{x_1} \edge[dd]_\pi &9\toplabel{x_9} \edge[d]_\mu
&10\toplabel{x_{10}} \edge[d]_\nu &&11\toplabel{x_{11}} \edge[d]_\psi
&12\toplabel{x_{12}} \edge[d]_\chi\\
 &7\edge[d]_\xi &8\edge[d]_\eta \edge[dr]^\tau &&6\edge[d]_\rho
&6\edge[d]_\rho\\ 2\threepool &2 &2 &5\threepool &5 &5 }$$ }

\noindent is to encode the following information: The module $A$ is generated
by top elements $x_i$ of type $e_i$ $(i=1,9,10,11,12)$ which are $K$-linearly
independent modulo $JA$ (here automatically satisfied, since $e_i\ne e_j$ for
$i\ne j$) such that $J^3A=0$, and 

(a) $JA/J^2A\cong S_7\oplus S_8\oplus S_6^2$, with $\mu x_9$, $\nu x_{10}$,
$\psi x_{11}$, $\chi x_{12}$ being $K$-linearly independent modulo $J^2A$;

(b) $J^2A\cong S_2^2\oplus S_5^2$, and the ``pooled elements'' $\pi x_1$,
$\xi\mu x_9$,
$\eta\nu x_{10}$ are
$K$-linearly dependent, while any two of these elements are $K$-linearly
independent; analogously, $\tau\nu x_{10}$, $\rho\psi x_{11}$, $\rho\chi
x_{12}$ are $K$-linearly dependent, with any subset of two $K$-linearly
independent.

Note that in this particular example, the module $A$ is determined up to
isomorphism by its graph, since the coefficients arising in the mentioned
linear dependence relations can be adjusted by suitably modifying the top
elements by scalar factors.

It is clear that we will not lose information if we omit the labels on edges
\ignore{$$\xymatrixrowsep{1pc}\xymatrix{ i\edge[d]\\ j}$$}  with the property
that there is a unique arrow $i\rightarrow j$ in $\Gamma$. Moreover, it
should be self-explanatory that certain countably generated left
$\la$-modules can be communicated by graphs as well. The {\it infinite\/}
graph

\ignore{
$$\xymatrixcolsep{0.8pc}\xymatrixrowsep{1.5pc}
\xymatrix{ 6\edge[dr]_(0.55)\sigma &&6\edge[dl]_(0.45)\rho
\edge[dr]_(0.55)\sigma &&6\edge[dl]_(0.45)\rho \edge[dr]_(0.55)\sigma
&&\cdots\cdots\\
 &5 &&5 &&5 &\cdots\cdots }$$  }

\noindent for instance, goes with a module $B$, uniquely determined up to
isomorphism, which is generated by top elements $x_1,x_2,x_3,\dots$ of type
$e_6$ which are $K$-linearly independent modulo $JB$ such that $\rho x_1=0$ 
and $\sigma x_i=
\rho x_{i+1}$ for all $i\ge1$, and such that the elements $\sigma x_i$,
$i\in\NN$, are
$K$-linearly independent (modulo $J^2B=0$).

\head 2. Contravariant finiteness and first examples\endhead

In 1980, Auslander and Smal\o\ \cite{\AuSm} introduced the following
definitions -- along with the duals -- in connection with their search for
conditions ensuring existence of almost split sequences. It turned out that
existence is guaranteed in full, extension-closed subcategories
$\A$ of
$\lamod$ which are both co- and contravariantly finite in the sense recalled
below. Even though their results have been  applied on numerous occasions
during the 1980's, it was not until the 1990's that the concept of
contravariant finiteness reached a high level of popularity, due to its links
with homology and tilting. The spark in the tinder barrel was a paper by
Auslander and Reiten
\cite{\AuRe}. We will describe their homological results in Section 3.

\definition{Definitions} Let $\A\subseteq\lamod$ be a full subcategory and
$X\in\lamod$. 

(1) A {\it right $\A$-ap\-prox\-i\-ma\-tion of $X$\/} is a homomorphism $f :
A\rightarrow X$ with $A\in\A$ such that each $g\in \Hom(B,X)$ with $B\in\A$
factors through
$f$, or equivalently, such that the following sequence of functors induced by
$f$ is exact:
$$\Hom_\la(-,A)|_\A \longrightarrow \Hom_\la(-,X)|_\A \longrightarrow 0.$$
(Since we will hardly mention the dual concept of `left approximation', we
will systematically suppress the qualifier `right' when discussing
approximations.)

(2) The subcategory $\A$ is said to be {\it contravariantly finite\/} (in
$\lamod$) if each $X\in\lamod$ has an $\A$-ap\-prox\-i\-ma\-tion, i.e., if
each of the
 functors $\Hom_\la(-,X)|_\A$ is finitely generated in the category of all
contravariant functors from $\A$ to {\bf Ab}.
\enddefinition

Suppose that $X$ has an $\A$-ap\-prox\-i\-ma\-tion. By a slight abuse of
language, we will then also refer to the source of this map as an
approximation. Not surprisingly, the $\A$-ap\-prox\-i\-ma\-tions of $X$ of
minimal length are all isomorphic. Indeed, as  was shown by Auslander and
Smal\o\
\cite{\AuSm}, given a minimal $\A$-ap\-prox\-i\-ma\-tion $f : A\rightarrow X$
and any
$\A$-ap\-prox\-i\-ma\-tion $f' : A'\rightarrow X$, there exists a split
embedding $g : A\rightarrow A'$ which makes the following diagram commute:

\ignore{
$$\xymatrixrowsep{1.5pc}\xymatrixcolsep{1pc} \xymatrix{ A\ar[rr]^g \ar[dr]_f
&&A'\ar[dl]^{f'}\\
 &X }$$ }

\noindent It is thus justified to refer to \underbar{the} minimal
$\A$-ap\-prox\-i\-ma\-tion of $X$ in case of existence.

If $\A$ is a resolving subcategory of $\lamod$, i.e., if $\A$ contains all
projectives in $\lamod$ and is closed under extensions and kernels of
epimorphisms, then the simples play a prominent role in checking
contravariant finiteness. Indeed:

\proclaim{Theorem 1} {\rm \cite{\AuRe}} Suppose $\A$ is a resolving
subcategory of $\lamod$. Then $\A$ is contravariantly finite in $\lamod$ if
and only if each of the simple left $\la$-modules has an
$\A$-ap\-prox\-i\-ma\-tion.
\qed\endproclaim

Since, clearly, our favorite category $\pinflamod$ is resolving, this will
provide a convenient test for contravariant finiteness. As we will see in the
next section, the minimal \pinfapprox s of the simples are the basic
structural building blocks for the objects in $\pinf(\Lamod)$ in case of
existence, which makes it a priority to understand the structure of these
particular approximations.

To provide  examples, we  begin with some

\subhead Well Known Facts\endsubhead

$\bullet$ A subcategory $\A\subseteq\lamod$ is contravariantly finite in case
it is `very big' or `very small'. Indeed, if $\A=\lamod$, then clearly
contravariant finiteness is guaranteed, the minimal approximations being the
identity maps. So, in particular, $\pinflamod$ is contravariantly finite
provided that $\la$ has finite global dimension. On the other hand, if
$\A$ has finite representation type, i.e., if there exist objects
$A_1,\dots,A_n\in\A$ such that each object in $\A$ is a direct sum of copies
of the $A_i$, we have contravariant finiteness as well \cite{\AuRe}. To
construct an
$\A$-ap\-prox\-i\-ma\-tion of a module $X$, simply add up as many copies of
each $A_i$ as the $K$-dimension of
$\Hom_\la(A_i,X)$ indicates. In particular, if $\A$ is the category of all
finitely generated projectives in $\lamod$, the minimal
$\A$-ap\-prox\-i\-ma\-tions are precisely the projective covers. This latter
category coincides with
$\pinflamod$ precisely when $\lfindim\la=0$.

$\bullet$ \cite{\AuRe} If $\la$ is stably equivalent to a hereditary algebra
$\la'$ (meaning that the stable category $\la\text{-\underbar{mod}}$, obtained
as a factor category from
$\lamod$ by killing the maps that factor through projectives, is equivalent to
the corresponding stable category  $\la'\text{-\underbar{mod}}$), then
$\pinflamod$ is contravariantly finite. This hypothesis is, in particular,
satisfied if $J^2=0$, and in that case, the minimal approximation $A(X)$ of
any
$\la$-module $X$ can be readily pinned down as follows \cite{\convenient}:
$A(X)= P/(JP)_{\text{fin}}$, where $P$ is the projective cover of $X$ and
$(JP)_{\text{fin}}$ is the direct sum of those homogeneous components of $JP$
which have finite projective dimension.

$\bullet$ \cite{\BuHZ} If $\la$ is left serial, then $\pinflamod$ is always
contravariantly finite. The minimal approximations of the simple left modules
arising in this situation will be described in the next section.

$\bullet$ The first example for which $\pinflamod$ was shown not to be
contravariantly finite is due to Igusa-Smal\o-Todorov
\cite{\IgSmTo}. It is a monomial relation algebra with $J^3=0$ and
$\lfindim\la=1$ which is closely related to the Kronecker algebra, its
$K$-dimension exceeding that of the latter algebra only by 2. In this example,
the right finitistic dimension is 0, which, in view of our first remark,
demonstrates that the right-hand category $\pinf(\text{mod-}\la)$
\underbar{is} contravariantly finite in $\text{mod-}\la$. Thus contravariant
finiteness of
$\pinf(-)$ is not left-right symmetric.

$\bullet$ In \cite{\HaHZ}, very general criteria for the failure of
contravariant finiteness of $\pinflamod$ were developed. We refer the reader
to the examples given there and in \cite{\HZSmdnsn}.

\medskip

The emerging picture indicates that, while hardly any of the traditionally
considered classes of finite dimensional algebras enjoy -- en bloc -- the
property that $\pinf(-)$ is contravariantly finite, the positive case is
ubiquitous. In fact, the condition of having contravariantly finite $\pinf(-)$
appears to slice diagonally through the prominent classes of finite
dimensional algebras.

\subhead First Installment of Nonstandard Examples\endsubhead

We will next present a first set of examples which are to informally
communicate prototypical phenomena ensuring that a given simple module  has a
\pinfapprox\ or fails to have such an approximation. These specific algebras
will then continue to serve us as illustrations along the way. We will follow
with proofs of some of the positive instances, but defer the discussion of
the negative instances to Section 4.

\definition{Examples A} Our first example shows that, for each natural number
$n$, there exists a finite dimensional monomial relation algebra $\la$ and a
simple
$S\in\lamod$ such that $\pinflamod$ is contravariantly finite and the minimal
\pinfapprox\ of $S$ is a direct sum of $n$ distinct nonzero indecomposable
components.
\enddefinition

\definition{{\bf A.1}} Fix $n\in\NN$, and let $\la=\kgami$ be the monomial
relation algebra with quiver $\Gamma$

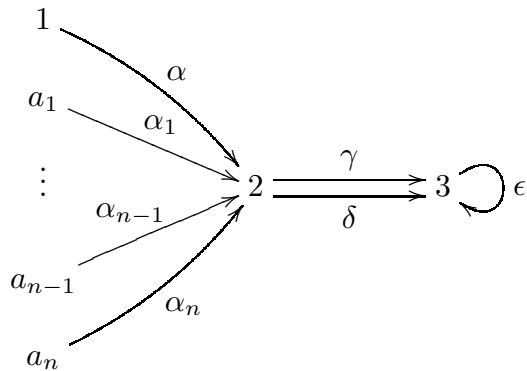
\ignore{
$$\xymatrixcolsep{4pc}\xymatrixrowsep{1pc}
\xymatrix{ 1\ar@/^/[ddr]^\alpha\\ a_1\ar[dr]^{\alpha_1}\\
\vdots &2\ar[r]<0.6ex>^\gamma \ar[r]<-0.6ex>_\delta &3\ar@(ur,dr)^\epsilon\\
a_{n-1}\ar[ur]^{\alpha_{n-1}}\\ a_n\ar@/_/[uur]_{\alpha_n} }$$ }

\noindent having indecomposable projective left modules with graphs

\ignore{
$$\xymatrixcolsep{0.8pc}\xymatrixrowsep{1.5pc}
\xymatrix{ 1\edge[d]_\alpha &&&2\edge[dl]_\gamma \edge[dr]^\delta
&&&3\edge[d]^\epsilon &&a_1\edge[d]^{\alpha_1} &&&&a_n\edge[d]^{\alpha_n}\\
2\edge[d]_\gamma &&3 &&3 &&3 &&2\edge[d]^\delta &&\cdots\cdots
&&2\edge[d]^\delta\\ 3 &&&&&&&&3 &&&&3 }$$ }

\noindent Then $\pinflamod$ is contravariantly finite, and the minimal
\pinfapprox s $A_1$, $A_2$, $A_3$, $A_{a_1}$, \dots, $A_{a_n}$ of the simple
modules $S_1$, $S_2$, $S_3$, $S_{a_1}$, \dots, $S_{a_n}$ have the following
graphs (which determine the corresponding modules up to isomorphism).

\ignore{
$$\xymatrixcolsep{0.5pc}\xymatrixrowsep{0.75pc}
\xymatrix{
 &1\edge[ddr] &&a_1\edge[ddl] &&&&1\edge[ddr] &&a_n\edge[ddl] &&&
&&2\edge[ddl]
\edge[ddr]\\ A_1: &&&&\bigoplus &\cdots &\bigoplus &&& &&&A_2:\\
 &&2 &&&&&&2 & &&& &3 &&3\\
 && &3\edge[dd] &&&& &a_i\edge[ddr] &&1\edge[ddl]\\
 &&A_3: & &&&&A_{a_i}: &&& &\save+<6ex,0ex> \drop{(1\le i\le n)}\restore\\
 && &3 &&&& &&2 }$$ }
\enddefinition

\definition{{\bf A.2}} Now let $\la= K\Gamma'/I'$, where $\Gamma'$ is obtained
from the quiver $\Gamma$ of A.1 by removing the arrow $\gamma$, and the ideal
$I'\subseteq K\Gamma'$ is such that the graphs of $\la e_3$, $\la e_{a_1}$,
\dots, $\la e_{a_n}$ are as under A.1, whereas $\la e_1$ and $\la e_2$ have
graphs

\ignore{
$$\xymatrixcolsep{1pc}\xymatrixrowsep{0.75pc}
\xymatrix{ 1\edge[dd]_\alpha &&&&2\edge[dd]^\delta\\
 &&\txt{and}\\ 2 &&&&3 }$$ }

\noindent respectively.

Then $S_1=\la e_1/Je_1$ has minimal \pinfapprox

\ignore{
$$\xymatrixcolsep{1pc}\xymatrixrowsep{1.5pc}
\xymatrix{ 1\edge[drr] &a_1\edge[dr] &a_2\edge[d] &\cdots\cdots
&a_n\edge[dll]\\
 &&2 }$$ }

\noindent $S_2$ has minimal \pinfapprox\ $\la e_2$, and
$S_{a_1}$, \dots, $S_{a_n}$ belong to $\pinflamod$. In particular,
$\pinflamod$ is again contravariantly finite.
\enddefinition

\definition{Examples B}
\enddefinition

\definition{{\bf B.1}} Let $\la=\kgami$, where $\Gamma$ is the quiver

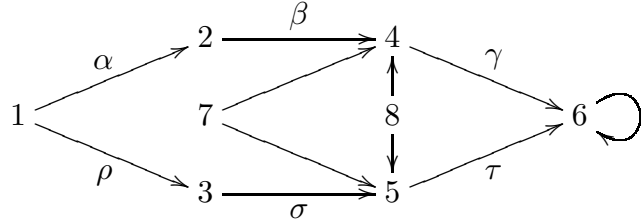
\ignore{
$$\xymatrixcolsep{4pc}\xymatrixrowsep{1pc}
\xymatrix{
 &2\ar[r]^\beta &4\ar[rd]^\gamma\\ 1\ar[ur]^\alpha \ar[dr]_\rho &7\ar[ur]
\ar[dr] &8\ar[u] \ar[d] &6\ar@(ur,dr)\\
 &3\ar[r]_\sigma &5\ar[ur]_\tau }$$ }

\noindent and the ideal $I\subseteq K\Gamma$ contains the relation
$\gamma\beta\alpha -\tau\sigma\rho$, together with suitable monomial
relations, such that the indecomposable projective left $\la$-modules have
graphs

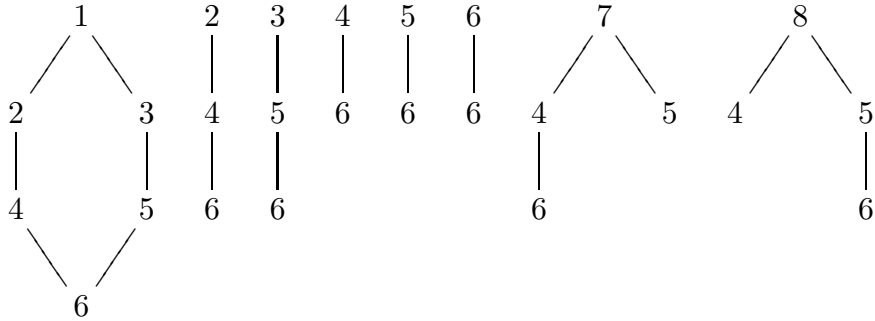
\ignore{
$$\xymatrixcolsep{0.8pc}\xymatrixrowsep{1.5pc}
\xymatrix{
 &1\edge[dl] \edge[dr] &&2\edge[d] &3\edge[d] &4\edge[d] &5\edge[d] &6\edge[d]
&&7\edge[dl] \edge[dr] &&&8\edge[dl] \edge[dr]\\ 2\edge[d] &&3\edge[d]
&4\edge[d] &5\edge[d] &6 &6 &6 &4\edge[d] &&5 &4 &&5\edge[d]\\ 4\edge[dr]
&&5\edge[dl] &6 &6 &&&&6 &&&&&6\\
 &6 }$$ }

\noindent Then $S_i\in \pinflamod$ for $i=2,3$, the minimal
\pinfapprox s of $S_4$, $S_5$, $S_6$ are $\la e_4$, $\la e_5$, and $\la e_6$,
respectively, while none of $S_1$, $S_7$, $S_8$ has a
\pinfapprox. In fact, in Section 4, we will see that there is no object
$A\in\pinflamod$ such that all of the homomorphisms from the modules in the
following subclass of $\pinflamod$ -- call it $\C$ -- to $S_1$ factor through
$A$.

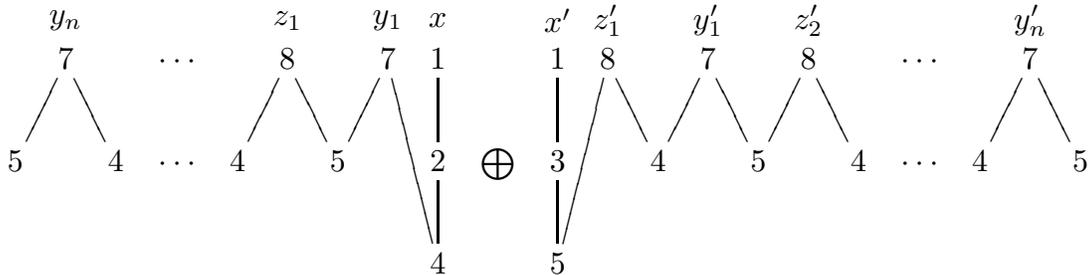
\ignore{
$$\xymatrixcolsep{0.4pc}\xymatrixrowsep{1.5pc}
\xymatrix{
 &7\toplabel{y_n} \edge[dl] \edge[dr] &&\cdots &&8\toplabel{z_1} \edge[dl]
\edge[dr] &&7\toplabel{y_1} \edge[dl]
\edge[ddr] &1\toplabel{x} \edge[d] &&1\toplabel{x'} \edge[d] &8\toplabel{z'_1}
\edge[ddl]
\edge[dr] &&7\toplabel{y'_1} \edge[dl] \edge[dr] &&8\toplabel{z'_2} \edge[dl]
\edge[dr] &&\cdots &&7\toplabel{y'_n} \edge[dl]
\edge[dr]\\  5 &&4 &\cdots &4 &&5 &&2\edge[d] &\bigoplus &3\edge[d] &&4 &&5
&&4 &\cdots &4 &&5\\
 &&&& &&&&4 &&5 }$$ }
\enddefinition

\definition{{\bf B.2}} Now let $\la$ be the factor algebra of the algebra
described in B.1, modulo the ideal generated by $e_7$ and $e_8$. Then the
graphs of the indecomposable projective modules $\la e_1$, \dots, $\la e_6$
remain unchanged, but this time $\pinflamod$ is contravariantly finite.
Indeed, the minimal
\pinfapprox s of $S_2$, \dots, $S_6$ are as above, while $S_1$ has the
following minimal \pinfapprox:

\ignore{
$$\xymatrixcolsep{1.0pc}\xymatrixrowsep{1.5pc}
\xymatrix{ 1\edge[d] &&1\edge[d]\\ 2\edge[d] &\bigoplus &3\edge[d]\\ 4 &&5 }$$
}
\enddefinition

The next example shows that the structure of the minimal
\pinfapprox s of the simple modules need by no means be as simplistic as in
the previous instances, not even in situations where the indecomposable
projective modules are of a simplistic makeup.

\definition{Examples C}
\enddefinition

\definition{{\bf C.1}} This time, let $\la=\kgami$ be the monomial relation
algebra with quiver $\Gamma$

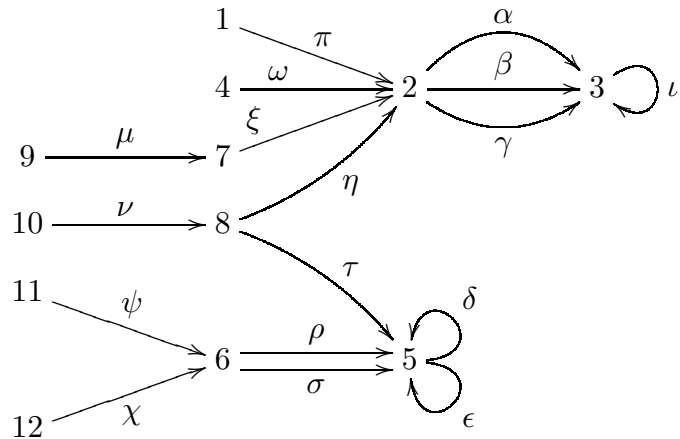
\ignore{
$$\xymatrixcolsep{4pc}\xymatrixrowsep{0.75pc}
\xymatrix{
 &1\ar[dr]^\pi\\
 &4\ar[r]^(0.3)\omega &2\ar[r]^\beta \ar@/^1.5pc/[r]^\alpha
\ar@/_1.0pc/[r]_\gamma &3\ar@(ur,dr)^\iota\\  9\ar[r]^\mu &7\ar[ur]^(0.2)\xi\\
10\ar[r]^\nu &8\ar@/_/[uur]_(0.55)\eta \ar@/^/[ddr]^(0.55)\tau\\
11\ar[dr]^\psi\\
 &6\ar[r]<0.6ex>^\rho \ar[r]<-0.6ex>_\sigma &5\ar@(r,u)_\delta
\ar@(r,d)^\epsilon\\ 12\ar[ur]_\chi }$$ }

\noindent and choose the ideal $I\subseteq K\Gamma$ of relations so that the
indecomposable projective left $\la$-modules have the following graphs:

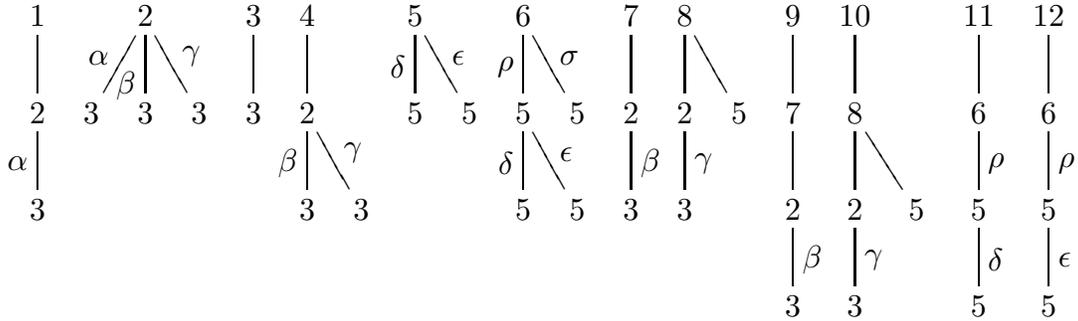
\ignore{
$$\xymatrixcolsep{0.5pc}\xymatrixrowsep{1.5pc}
\xymatrix{ 1\edge[d] &&2\edge[dl]_\alpha \edge[d]_(0.7)\beta \edge[dr]^\gamma
&&3\edge[d] &4\edge[d] &&5\edge[d]_\delta \edge[dr]^\epsilon &&6\edge[d]_\rho
\edge[dr]^\sigma &&7\edge[d] &8\edge[d] \edge[dr] &&9\edge[d] &10\edge[d]
&&11\edge[d] &12\edge[d]\\ 2\edge[d]_\alpha &3 &3 &3 &3 &2\edge[d]_\beta
\edge[dr]^\gamma &&5 &5 &5\edge[d]_\delta \edge[dr]^\epsilon &5
&2\edge[d]^\beta &2\edge[d]^\gamma &5 &7\edge[d] &8\edge[d] \edge[dr]
&&6\edge[d]^\rho &6\edge[d]^\rho\\ 3 &&&&&3 &3 &&&5 &5 &3 &3
&&2\edge[d]^\beta &2\edge[d]^\gamma &5 &5\edge[d]^\delta &5\edge[d]^\epsilon\\
 &&&&& & &&& & & & &&3 &3 &&5 &5  }$$ }

\noindent In this example, $\pinflamod$ is again contravariantly finite, and
the minimal approximations $A_1$, \dots, $A_{12}$ of the simples $S_1$, \dots,
$S_{12}$ are determined by their graphs as follows:

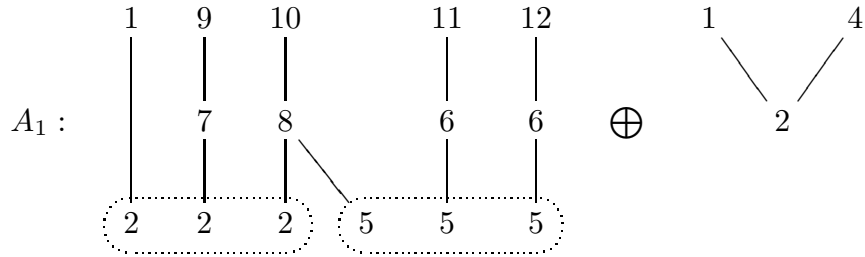
\ignore{
$$\xymatrixcolsep{1.0pc}\xymatrixrowsep{1.5pc}
\xymatrix{
 &1\edge[dd] &9\edge[d] &10\edge[d] &&11\edge[d] &12\edge[d] &&1\edge[dr]
&&4\edge[dl]\\ A_1: &&7\edge[d] &8\edge[d] \edge[dr] &&6\edge[d] &6\edge[d]
&\bigoplus &&2\\
 &2 \threepool &2 &2 &5 \threepool &5 &5 }$$ }
\ignore{
$$A_2\cong \Lambda e_2 \qquad\qquad\qquad A_3\cong \Lambda e_3
\qquad\qquad\qquad A_4: 
\xymatrixcolsep{1.0pc}\xymatrixrowsep{1.5pc}
\vcenter
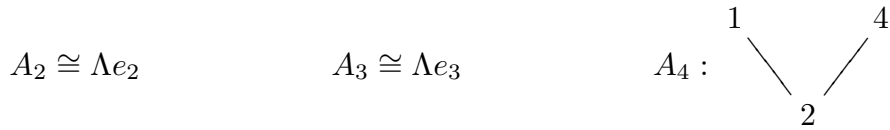
{\xymatrix{ 1 \edge[dr] &&4\edge[dl]\\
 &2 }}$$ }
\ignore{
$$\xymatrixcolsep{0.8pc}\xymatrixrowsep{1.5pc}
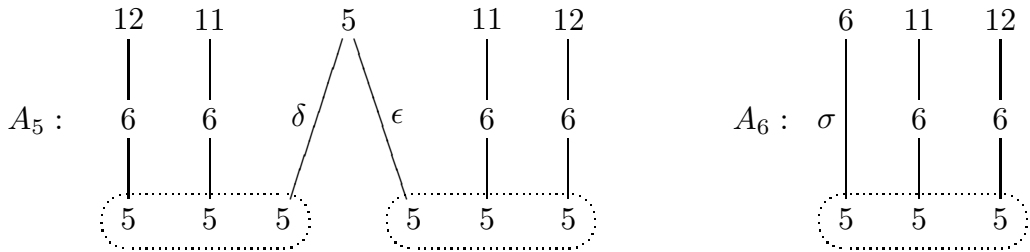
\xymatrix{
 &12\edge[d] &11\edge[d] &&5\edge[ddl]_\delta \edge[ddr]^\epsilon &&11\edge[d]
&12\edge[d] &&&&6\edge[dd]_\sigma &11\edge[d] &12\edge[d]\\ A_5: &6\edge[d]
&6\edge[d] &&&&6\edge[d] &6\edge[d] &&&A_6: &&6\edge[d] &6\edge[d]\\
 &5 \threepool &5 &5 &&5 \threepool &5 &5 &&&&5 \threepool &5 &5 }$$ }
\ignore{
$$\xymatrixcolsep{0.65pc}\xymatrixrowsep{1.5pc}
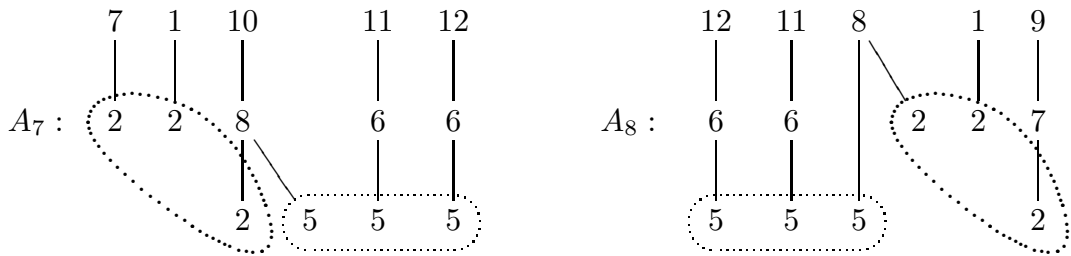
\xymatrix{
 &7\edge[d] &1\edge[d] &10\edge[d] &&11\edge[d] &12\edge[d] &&&&12\edge[d]
&11\edge[d] &8\edge[dd] \edge[dr] &&1\edge[d] &9\edge[d]\\ A_7:
 &2 \save[0,0]+(-3,0);[0,0]+(-3,0) **\crv{~*=<2.5pt>{.} [0,0]+(-3,4) 
&[0,1]+(4,4) &[1,2]+(4,0) &[1,2]+(3,-4) &[1,2]+(-3,-3) &[0,0]+(-3,-3)}\restore
 &2 &8\edge[d] \edge[dr] &&6\edge[d] &6\edge[d] &&&A_8: &6\edge[d] &6\edge[d] 
 &&2 \save[0,0]+(-3,0);[0,0]+(-3,0) **\crv{~*=<2.5pt>{.} [0,0]+(-3,4) 
&[0,1]+(4,4) &[1,2]+(4,0) &[1,2]+(3,-4) &[1,2]+(-3,-3) &[0,0]+(-3,-3)}\restore
 &2 &7\edge[d]\\
 &&&2 &5 \threepool &5 &5 &&&&5 \threepool &5 &5 &&&2  }$$ }
\ignore{
$$A_9\cong S_9 \qquad\qquad\qquad A_{10}\cong S_{10}
\qquad\qquad\qquad A_{11}\cong A_{12}: 
\xymatrixcolsep{1.0pc}\xymatrixrowsep{1.5pc}
\vcenter{\xymatrix{ 11\edge[d] &&12\edge[d]\\ 6\edge[dr] &&6\edge[dl]\\
 &5 }}$$ }
\enddefinition

\definition{{\bf C.2}} Finally, let $\la$ be the factor algebra of the algebra
under C.1, modulo the ideal generated by $e_{11}$ and $e_{12}$. Note that the
graphs of the indecomposable projectives $\la e_1$, \dots, $\la e_{10}$
remain the same as in C.1, since $e_{11}$ and $e_{12}$ are sources of
$\Gamma$. This time,
$\pinflamod$ fails to be contravariantly finite,
$S_1$,
$S_5$, $S_6$, $S_8$ being precisely those simples which have lost their
\pinfapprox s in the passage to the smaller algebra. As we will see in
Section 4, the homomorphisms onto $S_5$ from the modules of the following
$\pinflamod$-family  can, for instance, not all be factored through a fixed
$A\in\pinflamod$:

\ignore{
$$\xymatrixcolsep{0.5pc}\xymatrixrowsep{1.5pc}
\xymatrix{ 6\edge[dr]_(0.55)\sigma &&\cdots &&6\edge[dl]_(0.45)\rho
\edge[dr]_(0.55)\sigma &&6\edge[dl]_(0.45)\rho \edge[dr]_(0.55)\sigma
&&6\edge[dl]_(0.45)\rho
\edge[dr]_(0.55)\sigma &&6\edge[dl]_(0.45)\rho \edge[dr]_(0.55)\sigma &&\cdots
&&6\edge[dl]_(0.45)\rho \edge[dr]_(0.55)\sigma\\
 &5 &\cdots &5 &&5 &&5 &&5 &&5 &\cdots &5 &&5 }$$ }

\noindent (On the other hand, observe that they can be factored through the
module $A_6$ over the algebra in C.1.)
\enddefinition

To sketch justifications for some of the claims in the preceding examples, we
first spell out an obvious sufficient condition for a simple module
$S=\la e/Je$ to have a \pinfapprox. Indeed, this is the case provided that the
following is true: There exist indecomposable modules $T_1$,
\dots, $T_m$ in $\pinflamod$ with top elements $x_i\in T_i$ of type $e$ such
that for each indecomposable object $X\in\pinflamod$ having a top element $x$
of type $e$, there exists a factor module $X/Y$ with $\overline{x}= x+Y\ne 0$
which can be embedded into some $T_j$ in such a fashion that $\overline x$ is
mapped to $x_j$. If there exist $T_1$, \dots, $T_m$ as stipulated, we know
moreover that the indecomposable direct summands of the minimal
\pinfapprox s of $S$ are all recruited from the $T_i$.

\definition{Ad Example A.1} In arguing that $A_1$ is the minimal
\pinfapprox\ of $S_1$, we bypass the facts that $\lfindim\la=1$ (use the
method of \cite{\pre}, for instance), that the module $A_1$ has finite
projective dimension and satisfies the required minimality condition. To
verify that $A_1$ is a \pinfapprox\ of $S_1$, we let
$X\in\pinflamod$ be indecomposable and endowed with a top element $x$ of type
$e_1$. Then $\alpha x\ne0$, since otherwise the module with graph 
\ignore{$\vcenter{\xymatrixrowsep{1pc}\xymatrix{2\edge[d]\\ 3}}$} would be a
direct summand of $\Omega^1(X)$, which is impossible. If
$\gamma\alpha x\ne0$, then $X\cong \la x\cong \la e_1$ (use indecomposability
and finite projective dimension), and $X/\soc X$ embeds into each of the
modules 
\ignore{$\vcenter{\xymatrixrowsep{1pc} \xymatrixcolsep{0.5pc}
\xymatrix{1\edge[dr] &&a_i\edge[dl]\\ &2}}$}
 for $i=1,\dots,n$. So suppose that $\gamma\alpha x=0$. In this case,
$J^2X=0$, the socle of $X$ being homogeneous of type $e_2$, and
$\Omega^1(X)\cong (\la e_2)^r$ for some $r$. One infers the existence of a
factor module $X/Y$ with graph

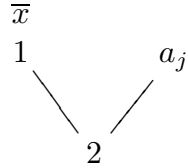
\ignore{
$$\xymatrixcolsep{1.0pc}\xymatrixrowsep{1.5pc}
\xymatrix{ 1\edge[dr] \save+<0ex,3ex> \drop{\overline x}\restore
&&a_j\edge[dl]\\
 &2 }$$ }

\noindent for some $j$. \qed\enddefinition

\definition{Ad A.2} To see that, also in this example, the exhibited module
$A_1$ is the minimal
\pinfapprox\ of $S_1$, observe that 
\ignore{$\Omega^1(A_1)= \biggl( \vcenter{\xymatrixrowsep{1pc}
\xymatrix{2\edge[d]^\delta\\ 3}} \biggr)^n =(\Lambda e_2)^n$},
 whence $A_1\in\pinflamod$. The rest of the argument is similar to the one
given above. \qed\enddefinition

\definition{Ad B.2} Once more, we will show that $A_1$ is a
\pinfapprox\ of $S_1$ (minimality being clear then). Any indecomposable
module in $\pinflamod$ containing $\la e_1$ is clearly isomorphic to $\la
e_1$, and the only proper nonzero factor modules of $\la e_1$ which embed
into indecomposable modules $X\in\pinflamod$ are the direct summands of
$A_1$, as well as

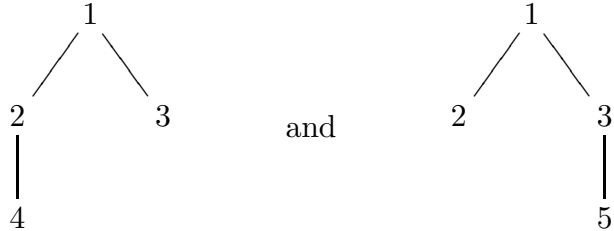
\ignore{
$$\xymatrixcolsep{1.0pc}\xymatrixrowsep{1.5pc}
\xymatrix{
 &1\edge[dl] \edge[dr] &&& &&&1 \edge[dl] \edge[dr]\\ 2\edge[d] &&3
&&\txt{and} &&2 &&3\edge[d]\\ 4 && && && &&5 }$$ }

\noindent Observe that none of these factors of $\la e_1$ has a  proper
extension to an indecomposable module in
$\pinflamod$. Since clearly the former factors are in turn factor modules of
the latter, our claim follows. \qed\enddefinition

All direct summands of the minimal approximations of the various simple
modules
$S=
\la e/Je$ exhibited above contain factor modules of $\la e$ which are minimal
 with respect to the property of being top-embeddable into modules of finite
projective dimension. (We say that a monomorphism
$A\rightarrow B$ is a {\it top-embedding\/} if it induces a monomorphism
$A/JA\rightarrow B/JB$.) This is generally true for the minimal \pinfapprox s
of simple modules, whenever $\la$ is either left serial or a string algebra
(see Sections 3 and 5), a fact which greatly facilitates resolving the
existence question (always positive in case of a left serial algebra and
algorithmically decidable for string algebras) and the construction of such
approximations. We conclude this section with an easy example showing that
this cannot be expected to hold in general. In Section 5, we will follow up
with an example demonstrating that the mentioned asset of string algebras is
not shared by arbitrary special biserial algebras either.

\definition{Example D} Let $\la=\kgami$ be the monomial relation algebra with
the following indecomposable projective left modules:

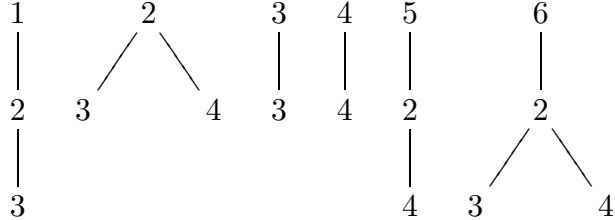
\ignore{
$$\xymatrixcolsep{0.8pc}\xymatrixrowsep{1.5pc}
\xymatrix{ 1\edge[d] &&2\edge[dl] \edge[dr] &&3\edge[d] &4\edge[d] &5\edge[d]
&&6\edge[d]\\ 2\edge[d] &3 &&4 &3 &4 &2\edge[d] &&2\edge[dl] \edge[dr]\\ 3
&&&&&&4 &3 &&4 }$$ }

\noindent Then $\pinflamod$ is contravariantly finite, the minimal
\pinfapprox\ $A_1$ of $S_1$ having the following graph:

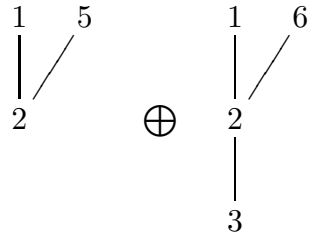
\ignore{
$$\xymatrixcolsep{0.8pc}\xymatrixrowsep{1.5pc}
\xymatrix{ 1\edge[d] &5\edge[dl] &&1\edge[d] &6\edge[dl]\\ 2 &&\bigoplus
&2\edge[d]\\
 &&&3 }$$ }

\noindent Note that the second summand of $A_1$ does not contain a proper
factor module of
$\la e_1$, whereas $\la e_1/\soc(\la e_1)$ is the (unique in this case)
minimal factor module of $\la e_1$ which can be top-embedded into an object of
$\pinflamod$. 
\enddefinition

\head 3. Homological importance of contravariant finiteness and a model
application\endhead

As surfaced in \cite{\dom}, \cite{\dep} and \cite{\Sma}, non-finitely
generated modules of finite projective dimension may display structural
phenomena which are completely different from those encountered in finitely
generated modules of finite projective dimension. (We labeled them `domino
effects' in
\cite{\dom}.) In particular, there may be objects in
$\pinf(\Lamod)$ whose projective dimension exceeds $\lfindim\la$ by any
predetermined positive number, even when $\lfindim\la=1$ \cite{\Sma}.
Moreover, the left cyclic finitistic dimension $\lcycfindim\la$, i.e., the
supremum of those projective dimensions which are attained on the cyclic
modules in
$\pinflamod$, may be strictly smaller than
$\lfindim\la$. In fact, for each natural number $n$, there exists a finite
dimensional algebra $\la_n$ such that $\lfindim\la_n$ is not attained on any
$n$-generated module (see \cite{\dep}). However, in case
$\pinflamod$ is contravariantly finite, all of the left finitistic dimensions
of
$\la$ coincide, and the objects of the big category
$\pinf(\Lamod)$ are as well understood as those of the small $\pinflamod$.

The following notation will be convenient: Given objects $A_1,\dots,A_n$ in
$\lamod$, let $\filt(A_1,\dots,A_n)$ be the full subcategory of $\lamod$ the
objects of which are those modules which have filtrations with consecutive
factors among $A_1,\dots,A_n$. More precisely, $X$ belongs to
$\filt(A_1,\dots,A_n)$ if and only if there exists a chain $X=X_0\supseteq
X_1\supseteq \cdots\supseteq X_m=0$ such that each of the factors
$X_i/X_{i+1}$ is isomorphic to some $A_{j(i)}$. Moreover,
$\Filt(A_1,\dots,A_n)$ will stand for the closure of $\filt(A_1,\dots,A_n)$
under direct limits in
$\Lamod$.

Concerning the structure of the {\it finitely generated\/} modules of finite
projective dimension in case $\A= \pinflamod$ is contravariantly finite,
Auslander and Reiten proved the following result in the more general context
of an arbitrary resolving subcategory $\A$.

\proclaim{Theorem 2} {\rm \cite{\AuRe}} Suppose that $\A$ is a resolving
contravariantly finite subcategory of $\lamod$, and that $A_1,\dots,A_n$ are
the minimal $\A$-ap\-prox\-i\-ma\-tions of the simple left $\la$-modules.
Then a module $X$ belongs to $\A$ if and only if $X$ is a direct summand of
an object in $\filt(A_1,\dots,A_n)$. \qquad$\square$\endproclaim

Of course, this theorem, applied to $\A=\pinflamod$, yields the following
consequence:

\proclaim{Corollary 3} {\rm \cite{\AuRe}} If $\pinflamod$ is contravariantly
finite and $A_1,\dots,A_n$ are the minimal
$\pinflamod$-ap\-prox\-i\-ma\-tions of the simple left $\la$-modules, then 
$$\lfindim\la= \lcycfindim\la= \max\{\pdim A_1,\dots,\pdim A_n\}.\qed$$
\endproclaim

Due to the fact that contravariant finiteness of $\pinflamod$ places demands
only on the finitely generated modules, the strong impact which this condition
has on {\it non-finitely generated\/} modules may come as a surprise. In fact,
the structure theory for objects in
$\pinflamod$ extends smoothly to
$\pinf(\Lamod)$.

\proclaim{Theorem 4} {\rm \cite{\HZSm}} Again suppose that $\pinflamod$ is
contravariantly finite, and let $A_1,\dots,\allowmathbreak A_n$ be as in the corollary. Then 
$$\pinf(\Lamod)= \Filt(A_1,\dots,A_n),$$ and, in particular,
$$\lFindim\la= \lfindim\la= \max\{\pdim A_1,\dots,\pdim A_n\}.\qed$$
\endproclaim

Thus, contravariant finiteness of $\pinflamod$ resolves the notorious
quandary of locating objects in $\Lamod$ on which $\lFindim\la$ is attained;
this search may be a very difficult task, even when finiteness of
$\lFindim\la$ is guaranteed in advance. The helpfulness of the above theory
will be displayed to full advantage in our examples.

\subhead Examples of Section 2 revisited\endsubhead

For the moment, we will only determine the finitistic dimensions of those
algebras displayed which give rise to contravariantly finite categories
$\pinf(-)$. By the preceding discussion, the objects of $\Lamod$ having finite
projective dimension are precisely the direct limits of the objects in
$\filt(A_1,\dots,A_n)$ with the $A_i$ as shown in Section 2.

\definition{Ad A.1} It is clear that $\pdim A_1= \pdim A_{a_i}=1$ for
$i=1,\dots,n$, whereas $\pdim A_2= \pdim A_3 =0$. Hence,
$\lFindim\la=\lfindim\la=1$.
\enddefinition

\definition{Ad A.2} In this example, the minimal \pinfapprox\ of each simple
$S_{a_i}$ centered in the vertex $a_i$  coincides with 
$S_{a_i}$ and has projective dimension 1, as does the minimal approximation of
$S_1$. The minimal approximations of
$S_2$ and $S_3$ are again identical with their projective covers. So, once
more, $\lFindim\la= \lfindim\la=1$.
\enddefinition

\definition{Ad B.2} Here $A_1$, $A_4$, $A_5$, $A_6$ are projective, while
$\pdim A_2= \pdim A_3=1$, and we obtain the same conclusion as before.
\qed\enddefinition

We will briefly digress from our main line of thought for another corollary of
the preceding theorem. The existence theorem for internal almost split
sequences \cite{\AuSm} which we quoted at the outset can be strengthened for
the category $\A=\pinflamod$ as follows.

\proclaim{Corollary 5} If $\pinflamod$ is contravariantly finite, then
$\pinflamod$ is also covariantly finite and thus has almost split sequences.
\endproclaim

\demo{Proof} By a result of Crawley-Boevey \cite{\CB}, it suffices to show
that arbitrary direct products of objects in the category $\pinflamod$ belong
to its closure $\overarrow{\pinf}(\lamod)$ under direct limits in
$\Lamod$. But in view of the theorem, contravariant finiteness of $\pinflamod$
entails the equality
$\pinf(\Lamod)=
\overarrow{\pinf}(\lamod)$ and, in view of the finiteness of $\lFindim\la$,
this guarantees closedness of this latter subcategory under direct products.
\qed\enddemo

As a class of examples of algebras $\la$ with very rich and complex module
categories, for which the above program of zeroing in on the structure of the
objects in
$\pinf(\Lamod)$ works to perfection, we will present the left serial algebras.
Recall that a split algebra $\la=\kgami$ is called {\it left serial\/} in case
no more than one arrow leaves any given vertex of $\Gamma$; equivalently,
this means that the indecomposable projective left $\la$-modules are all
uniserial. To describe the minimal \pinfapprox s of the simple modules in
this situation, we require the following definition.

\definition{Definition} Suppose that $T_1, \ldots , T_m$ is a sequence of
nonzero uniserial left $\la$-modules, and let $p_i$ be a {\it mast\/} of
$T_i$, namely a path in $K\Gamma$ of maximal length with $p_iT_i\ne0$. A left
$\la$-module
$T$ is called a  {\it saguaro on $(p_1, \ldots , p_m)$\/} if 

(i) $T\cong \bigl( \bigoplus _{1\leq i\leq m}T_i \bigr) \big/U$, where
$U\subseteq \bigoplus _{1\leq i\leq m}JT_i$ is generated by a sequence of
elements of  the form
$q_it_i-q'_{i+1}t_{i+1}$, $1\leq i\leq m-1$, where $t_i \in T_i$ are suitable
top elements and  $q_i, q'_i$ are right subpaths of the masts $p_i$ such that
$q_it_i\neq 0$, and
$q'_{i+1}t_{i+1}\neq 0$;  moreover, we require that 

(ii) each $T_j$ embeds canonically into $T$ via
$$T_j@>{\roman{can}}>> \biggl( \bigoplus_{1\leq i\leq m}T_i\biggr) \bigg/
U\cong T.$$ The uniserial modules $T_i$ are called the {\it trunks\/} of
$T$.

We will identify $T$ with
$\bigl( \bigoplus _{1\leq i\leq m}T_i \bigr) \big/ U$. To avoid ambiguities,
we will denote the canonical images of the trunks $T_i$ inside $T$ by $\That
_i$ and the canonical images of the top elements $t_i$ by $\that _i$.  Any
such sequence $(\that_1,
\ldots , \that_m)$ will be called a {\it canonical sequence of top elements\/}
for $T$.
\enddefinition

Note that saguaros are particularly amenable to graphing, the shape of their
graphs explaining their name (they share shape and name with a cactus found
in the Sonoran desert, {\it Cereus giganteus\/}). In fact, the definition
forces them to be glued together in a very straightforward fashion from their
uniserial trunks: Layered and labeled graphs relative to a canonical sequence
of top elements always exist (not only over left serial algebras), and are
built on the pattern illustrated below.

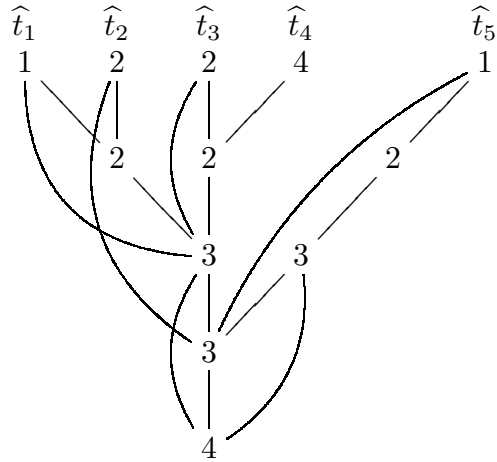
\ignore{
$$\xymatrixcolsep{1.5pc}\xymatrixrowsep{1.5pc}
\xymatrix{ 1\toplabel{\that_1} \edge[dr] \edge@/_1.75pc/[ddrr]
&2\toplabel{\that_2}
\edge[d]
\edge@/_1.75pc/[dddr] &2\toplabel{\that_3}
\edge[d]
\edge@/_1pc/[dd] &4\toplabel{\that_4} \edge[dl] &&1\toplabel{\that_5}
\edge[dl]
\edge@/_1pc/[dddlll]\\
 &2 \edge[dr] &2 \edge[d] &&2 \edge[dl]\\
 &&3 \edge[d] \edge@/_1pc/[dd] &3 \edge[dl] \edge@/^1pc/[ddl]\\
 &&3 \edge[d]\\
 &&4}$$ }

\noindent Here $T= \bigl( \bigoplus_{i=1}^5 T_i \bigr) \big/U$, where the
trunks $T_i= \la t_i$ of $T$ have graphs

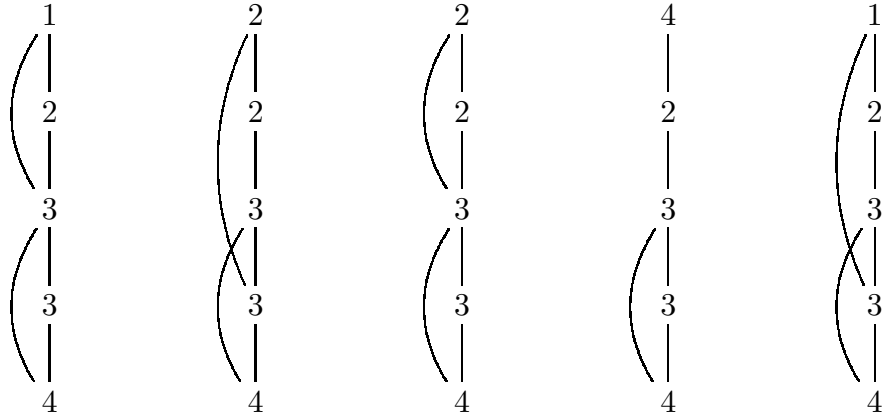
\ignore{
$$\xymatrixrowsep{1.5pc}
\xymatrix{ 1 \edge[d] \edge@/_1pc/[dd] &&2 \edge[d] \edge@/_1pc/[ddd] &&2
\edge[d]
\edge@/_1pc/[dd] &&4
\edge[d] &&1 \edge[d] \edge@/_1pc/[ddd]\\ 2 \edge[d] &&2 \edge[d] &&2
\edge[d] &&2
\edge[d] &&2 \edge[d]\\ 3 \edge[d] \edge@/_1pc/[dd] &&3 \edge[d]
\edge@/_1pc/[dd] &&3 \edge[d]
\edge@/_1pc/[dd] &&3 \edge[d] \edge@/_1pc/[dd] &&3 \edge[d] \edge@/_1pc/[dd]\\
3 \edge[d] &&3 \edge[d] &&3 \edge[d] &&3 \edge[d] &&3 \edge[d]\\
4&&4&&4&&4&&4}$$ }

\noindent relative to the top elements $t_i$.

Of course, over left serial algebras, the graphs of the uniserial left
modules are necessarily edge paths (in particular, uniserials have unique
masts), so that the graphs of saguaros simplify to trees. It turns out that,
in this situation, all minimal
\pinfapprox s of simple modules are recruited from the class of saguaros with
simple socles. In fact, these approximating saguaros can be pinned down
explicitly and even constructed algorithmically on the basis of quiver and
relations.

\proclaim{Theorem 6} {\rm \cite{\BuHZ}} Suppose that
$\la=\kgami$ is a left serial algebra. Then $\pinflamod$ is contravariantly
finite, and the minimal
\pinfapprox s of the simple left $\la$-modules are saguaros with simple
socles.

More precisely, the minimal
\pinfapprox\ of a simple left $\la$-module $S=\la e/Je$ can be described as
follows: If $\la e/C$ is the (unique) minimal nonzero factor module of $\la
e$ which has finite projective dimension, there is a unique saguaro $A(S)$ of
maximal length in
$\pinflamod$ such that $\la e/C$ is a trunk of $A(S)$ and $\soc A(S)$ is
simple. Moreover, the canonical epimorphisms $A(S)\rightarrow S$, which map
$\la e/C$ onto
$S$ and send the other trunks of $A(S)$ to zero, are minimal
\pinfapprox s. \qed\endproclaim

To refer back to the concluding remark of Section 2: In the setting of the
theorem, $\la e/C$  actually coincides with the nonzero factor module of $\la
e$ which is minimal with respect to top-embeddability into a module of finite
projective dimension.

Actually, not only is $\pinflamod$ always contravariantly finite in the left
serial case, but so are the categories $\Cal P^{(d)}= \Cal P^{(d)}(\lamod)$
consisting of the finitely generated left $\la$-modules of projective
dimensions at most $d$. Moreover, the minimal  $\Cal
P^{(d)}$-ap\-prox\-i\-ma\-tions of the simples are again saguaros, and the
sequences of these saguaros for $1\le d\le \lfindim\la$ record the homological
properties of $\la$ with high precision, the case $d=\lfindim\la$ leading back
to $\pinflamod$.

\definition{Example E}
Let $\la$ be a left serial algebra whose indecomposable projective modules
are represented by the following graphs.

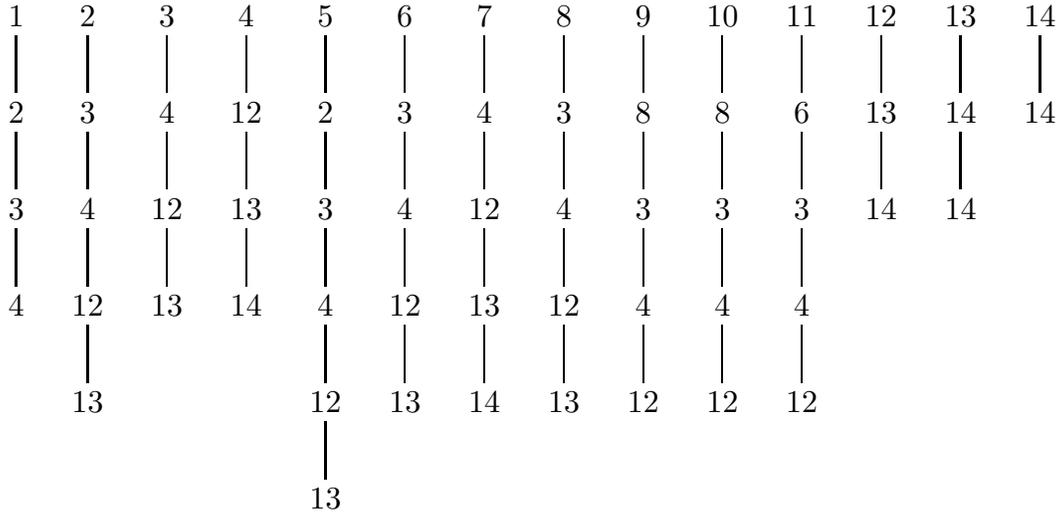
\ignore{
$$\xymatrixcolsep{0.75pc}\xymatrixrowsep{1.5pc}
\xymatrix{ 1 \edge[d] &2 \edge[d] &3 \edge[d] &4 \edge[d] &5 \edge[d] &6
\edge[d] &7
\edge[d] &8 \edge[d] &9 \edge[d] &10 \edge[d] &11 \edge[d] &12 \edge[d] &13
\edge[d] &14 \edge[d] \\ 2 \edge[d] &3 \edge[d] &4 \edge[d] &12 \edge[d] &2
\edge[d] &3 \edge[d] &4
\edge[d] &3 \edge[d] &8 \edge[d] &8 \edge[d] &6 \edge[d] &13 \edge[d] &14
\edge[d] &14\\ 3 \edge[d] &4 \edge[d] &12 \edge[d] &13 \edge[d] &3 \edge[d]
&4 \edge[d] &12
\edge[d] &4 \edge[d] &3 \edge[d] &3 \edge[d] &3 \edge[d] &14 &14\\ 4 &12
\edge[d] &13 &14 &4 \edge[d] &12 \edge[d] &13 \edge[d] &12 \edge[d] &4
\edge[d] &4 \edge[d] &4 \edge[d] \\
 &13 &&&12 \edge[d] &13 &14 &13 &12 &12 &12\\ &&&&13 }$$ }

\noindent The evolution of the ${\Cal P}^{(d)}$-ap\-prox\-i\-ma\-tions of the
simple left $\la$-module $S_1$ is graphically represented below.  We exhibit
the minimal 
${\Cal P}^{(1)}$-, ${\Cal P}^{(2)}$-, ${\Cal P}^{(3)}$-ap\-prox\-i\-ma\-tions
of $S_1$ from left to right; the last coincides with the minimal
\pinfapprox, since the left finitistic dimension of $\la$ is $3$. 

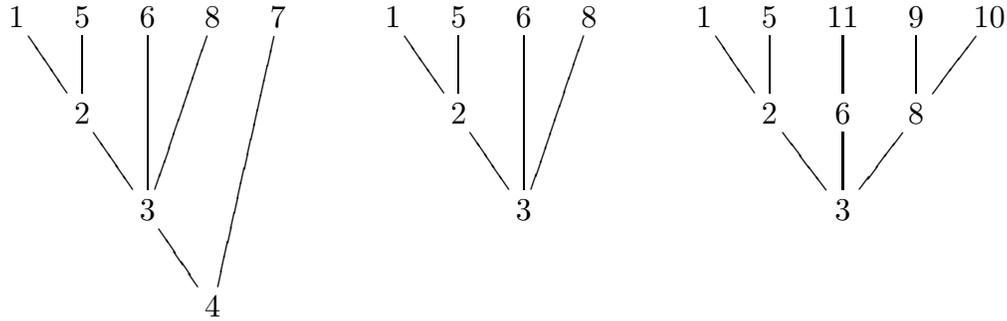
\ignore{
$$\xymatrixcolsep{0.8pc}\xymatrixrowsep{1.5pc}
\xymatrix{ 1 \edge[dr] &5 \edge[d] &6 \edge[dd] &8 \edge[ddl] &7 \edge[dddl]
&&1
\edge[dr] &5 \edge[d] &6 \edge[dd] &8 \edge[ddl] &&1 \edge[dr] &5 \edge[d]
&11 \edge[d] &9 \edge[d] &10 \edge[dl]\\
 &2 \edge[dr] &&&&&&2 \edge[dr] &&&&&2 \edge[dr] &6 \edge[d] &8
\edge[dl]\\
 &&3 \edge[dr] &&&&&&3 &&&&&3\\
 &&&4 }$$ }
\enddefinition

\head 4. Phantoms. Definitions, existence, and basic properties\endhead

The objects discussed in this section were introduced by Happel and the author
in \cite{\HaHZ}. As mentioned in the introduction, their purpose is twofold:
In the first place, they serve as indicators as to whether or not a given
subcategory $\A\subseteq\lamod$ is contravariantly finite. Their second role
is that of retaining the kind of information which is stored in minimal
$\A$-ap\-prox\-i\-ma\-tions whenever they exist, within potentially infinite
dimensional frames; this role is played most satisfactorily by the
`effective' phantoms.

Since the concept of a phantom is possibly not easily translated into an
intuitive picture, we break the definitions into several parts, and add a
fairly extensive discussion at each step. The subsequent detailed analysis of
the examples of Section 2 should also help the reader to see that we are
dealing with objects that arise very naturally when the relations of a module
$X$ are being compared with those of the objects in a subcategory of $\lamod$.

\definition{Definition, Part I} (Relative approximations) Let $\A$ be a full
subcategory of $\lamod$, and $\C$ a subcategory of $\A$. Moreover, let
$X\in\lamod$. 

A {\it $\C$-ap\-prox\-i\-ma\-tion of $X$ inside $\A$\/} is a homomorphism $f :
A\rightarrow X$ with $A\in\A$ such that each map in $\Hom_{\la}(C,X)$ with
$C\in\C$ factors through $f$. Again, we will loosely refer to the object $A$
as a
$\C$-ap\-prox\-i\-ma\-tion of $X$ inside $\A$.
\enddefinition

Clearly, whenever $X$ has a minimal $\A$-ap\-prox\-i\-ma\-tion, $A_0$ say, the
classical $\A$-ap\-prox\-i\-ma\-tions of $X$ are precisely the
$\{A_0\}$-ap\-prox\-i\-ma\-tions of $X$ inside $\A$. In particular, we then
obtain
$A_0$ as a $\C$-ap\-prox\-i\-ma\-tion of $X$ inside $\A$, where $\C$ is a
\underbar{finite} subcategory of $\A$. Moreover,
$A_0$ is a direct summand of any $\{A_0\}$-ap\-prox\-i\-ma\-tion of $X$
inside $\A$, and so, a fortiori, is a subfactor of any such approximation. On
the other hand, given a finite subcategory $\C\subseteq\A$, the module $X$
will have $\C$-ap\-prox\-i\-ma\-tions inside $\A$ provided that we require 
$\A$ to be closed under finite direct sums: Just sum up a sufficient number of
copies of the objects in $\C$. In other words, approximations relative to
\underbar{finite} classes are always available, also in case there are no
classical
$\A$-ap\-prox\-i\-ma\-tions of $X$.

\definition{Definition, Part II} (Phantoms) Retain the notation of Part I, and
suppose, in addition, that $\A$ is closed under finite direct sums.

A finitely generated module $H\in\lamod$ is an {\it $\A$-phan\-tom of $X$\/}
in case 

(*) there is a finite subcategory $\A'\subseteq\A$ such that $H$ occurs as a
subfactor of {\it every\/} $\A'$-ap\-prox\-i\-ma\-tion of $X$ inside $\A$. 

More generally, an arbitrary module $H\in\Lamod$ will be called an
$\A$-phan\-tom of
$X$ if each of its finitely generated submodules satisfies (*). Of course,
the choices of the finite subclasses  $\A'\subseteq\A$ will vary with the
finitely generated submodules $H'$ of $H$.
\enddefinition

The class of all $\A$-phan\-toms of $X$ is clearly closed under subfactors,
and consequently is closed under direct limits of direct systems as well. In
fact, a module $H\in\Lamod$ is an
$\A$-phan\-tom of
$X$ if and only if
$H$ is the direct limit of a direct system of finitely generated
$\A$-phan\-toms of $X$. This may make the class of $\A$-phan\-toms enormous:
Indeed, the class of all $\A$-phan\-toms of a simple module $S$ may encompass
the entire class of indecomposables $H$ in $\Lamod$ with $S\hookrightarrow
H/JH$. This slack in the definition of phantoms has the advantage of
facilitating their construction. Often, the particular structure of phantoms
is irrelevant -- their sheer size is enough to signal non-existence of
traditional
$\A$-ap\-prox\-i\-ma\-tions of $X$. In fact, if $X$ has an
$\A$-ap\-prox\-i\-ma\-tion then the class of $\A$-phan\-toms of $X$ coincides
with the set of subfactors of the minimal such approximation. Consequently,
the existence of phantoms of $X$ of unbounded lengths is enough to guarantee
non-existence of traditional
$\A$-ap\-prox\-i\-ma\-tions. Of course, in terms of encapsulating further
information, the usefulness of phantoms so generously defined is moderate.
Hence, we single out a subclass of phantoms which are more strongly tied to
the category $\A$ and carry a full complement of information on how a given
subcategory $\C \subseteq \A$ relates to $X$.  

\definition{Definition, Part III} (Effective phantoms) Keep the notation of
Part II, and denote by $\overarrow \A$ the closure of $\A$ under direct
limits (of direct systems) in
$\Lamod$. Moreover,  fix a subcategory
$\C\subseteq\A$.

An $\A$-phan\-tom $H\in\overarrow \A$ is called {\it effective relative to
$\C$\/} if there exists a homomorphism $h : H\rightarrow X$ with the property
that each map in $\Hom_{\la}(C,X)$ with $C\in\C$ factors through $h$. 
\enddefinition

In a self-explanatory extension of of the term `approximation' to
subcategories of $\Lamod$, the effectiveness condition thus calls for
$h : H \rightarrow X$ to be a $\C$-ap\-prox\-i\-ma\-tion of $X$ inside
$\overarrow \A$. In case $X$ has a traditional $\A$-ap\-prox\-i\-ma\-tion,
the minimal such approximation is clearly the only effective $\A$-phan\-tom
of $X$ relative to
$\A$. Otherwise, existence of interesting phantoms, effective or not, is not
immediately clear, but is guaranteed by the following result.

\proclaim{Theorem 7} {\rm \cite{\HaHZ}} Let $\A\subseteq\lamod$ be a full 
subcategory which is closed under finite direct sums. For $X\in\lamod$, the
following conditions are equivalent:

{\rm (1)} $X$ fails to have an $\A$-ap\-prox\-i\-ma\-tion.

{\rm (2)} $X$ has $\A$-phan\-toms of infinite $K$-dimension.

{\rm (3)} There exist  countable subclasses $\C\subseteq\A$ such that $X$ has
infinite dimensional $\A$-phan\-toms which are effective relative to $\C$.
\qed\endproclaim

Next we discuss our `test examples' in light of the new concepts.

\definition{Ad B.1} Here are two infinite dimensional $\pinflamod$-phan\-toms
of
$S_1$:

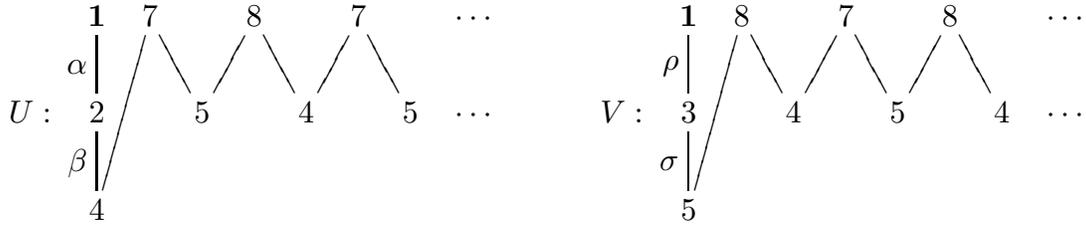
\ignore{
$$\xymatrixcolsep{0.45pc}\xymatrixrowsep{1.5pc}
\xymatrix{
 &{\bold 1}\edge[d]_\alpha &7\edge[ddl] \edge[dr] &&8\edge[dl] \edge[dr]
&&7\edge[dl]
\edge[dr] &&\cdots &&&&{\bold 1}\edge[d]_\rho &8\edge[ddl] \edge[dr]
&&7\edge[dl]
\edge[dr] &&8\edge[dl] \edge[dr] &&\cdots\\ U: &2\edge[d]_\beta &&5 &&4 &&5
&\cdots &&&V: &3\edge[d]_\sigma &&4 &&5 &&4 &\cdots\\
 &4 &&&&&&& &&&&5 }$$ }

\noindent Note that $U\oplus V$ is in turn a $\pinflamod$-phan\-tom of $S_1$,
since there is no object in $\pinf(\Lamod)$ having a graph with subgraph

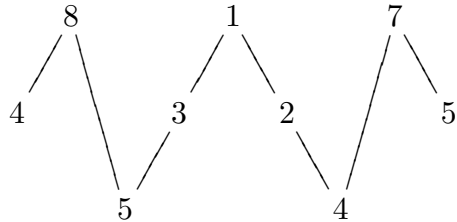
\ignore{
$$\xymatrixcolsep{0.5pc}\xymatrixrowsep{1.5pc}
\xymatrix{
 &8\edge[dl] \edge[ddr] &&&1\edge[dl] \edge[dr] &&&7\edge[ddl] \edge[dr]\\ 4
&&&3\edge[dl] &&2\edge[dr] &&&5\\
 &&5 &&&&4 }$$ }

\noindent In fact, the phantom $U\oplus V$ is effective relative to the class 
$\C\subseteq\pinflamod$ exhibited in B.1.

We will justify only that $U$ is a $\pinflamod$-phan\-tom of $S_1$. Consider
the class of objects

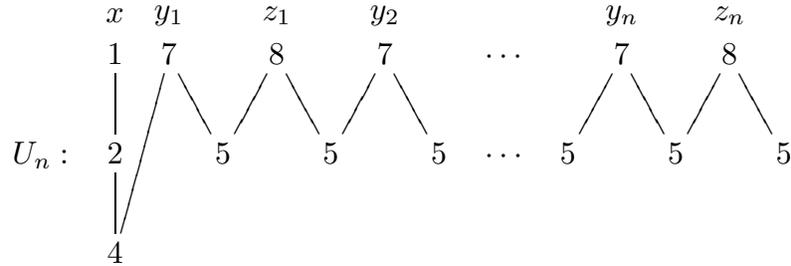
\ignore{
$$\xymatrixcolsep{0.5pc}\xymatrixrowsep{1.5pc}
\xymatrix{
 &1\toplabel{x} \edge[d] &7\toplabel{y_1} \edge[ddl] \edge[dr]
&&8\toplabel{z_1} \edge[dl] \edge[dr] &&7\toplabel{y_2} \edge[dl] \edge[dr]
&&\cdots &&7\toplabel{y_n} \edge[dl] \edge[dr] &&8 \toplabel{z_n} \edge[dl]
\edge[dr]\\   U_n: &2\edge[d] &&5 &&5 &&5 &\cdots &5 &&5 &&5\\
 &4 }$$ }

\noindent in $\pinflamod$, and observe that $U=\varinjlim U_n$. Hence, it
suffices to show that, for each $n\in\NN$, the module $U_n$ is a submodule of
each
$\{U_n\}$-ap\-prox\-i\-ma\-tion of $S_1$ inside $\pinflamod$. To see this,
fix $n$, and let $A$ be any $\{U_n\}$-ap\-prox\-i\-ma\-tion of $S_1$ inside
$\pinflamod$; say
$f : A\rightarrow S_1$ has the factorization property of the definition, and
$g\in
\Hom(U_n,A)$ factors the canonical epimorphism $U_n\rightarrow S_1$. Let
$a=g(x)$, where $x$ is the unique top element of type $e_1$ of $U_n$. Then
$a$ is a top element of
$A$. Due to the fact that finiteness of the projective dimension of $A$
entails either $\beta\alpha a\ne0$ or $\sigma\rho a\ne0$, our factorization
requirement forces
$\beta\alpha a$ to be nonzero. Consequently, $g(\beta\alpha x)= \beta\alpha
a$. If the arrows $7\rightarrow 4$ and $7\rightarrow 5$ are named $\chi$ and
$\psi$, respectively, we deduce $\beta\alpha a=f(\chi y_1)\ne0$, and hence
$f(y_1)\ne0$. The element $b_1=f(y_1)$ is a top element of $A$ (necessarily of
type $e_7$), since 7 is a source of $\Gamma$. To prevent the syzygy
$\Omega^1(A)$ from having a summand $S_5$ (the latter being incompatible with
finite projective dimension), we require $0\ne \psi(b_1)= \mu f(z_1)$, where
$\mu$ is the arrow $8\rightarrow 5$ and $z_1$ is as indicated in the above
figure. Setting
$c_1=f(z_1)$ and iterating this type of argument, we see that $0\ne \nu c_1=
f(\chi y_2)$, where
$\nu$ is the arrow
$8\rightarrow 4$. Set $b_2=f(y_2)$. It is readily checked that the top
elements $b_1,b_2$ are $K$-linearly independent modulo $JA$, and an obvious
induction on $m\le n$ gives us sequences of top elements, $b_1,\dots,b_n$ of
type
$e_7$, and
$c_1,\dots,c_n$ of type $e_8$ in $A$, both series being $K$-linearly
independent modulo $JA$. It is now straightforward to deduce that the
submodule of
$A$ generated by $a$ and the $b_i,c_i$ has the same graph as $U_n$. But this
graph clearly determines the corresponding module up to isomorphism, which
completes the argument.

Each module from the subclass $\D= \{D_n\mid n\in\NN\}$ of $\pinflamod$, with
$D_n$ determined by the graph

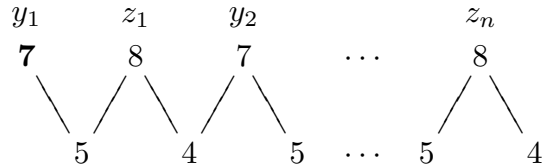
\ignore{
$$\xymatrixcolsep{0.5pc}\xymatrixrowsep{1.5pc}
\xymatrix{ {\bold 7}\toplabel{y_1} \edge[dr] &&8\toplabel{z_1} \edge[dl]
\edge[dr] &&7\toplabel{y_2} \edge[dl] \edge[dr] &&\cdots &&8\toplabel{z_n}
\edge[dl]
\edge[dr]\\
 &5 &&4 &&5 &\cdots &5 &&4 }$$ }

\noindent is a $\pinflamod$-phan\-tom of $S_7$, and consequently $\varinjlim
D_n$ provides us with an infinite dimensional $\pinflamod$-phan\-tom of $S_7$
which, moreover, is effective relative to the class $\D$. The simple module
$S_8$ behaves similarly.
\enddefinition

\definition{Ad C.2} An infinite dimensional $\pinflamod$-phan\-tom of $S_1$ is
determined (uni\-que\-ly, up to isomorphism) by the graph

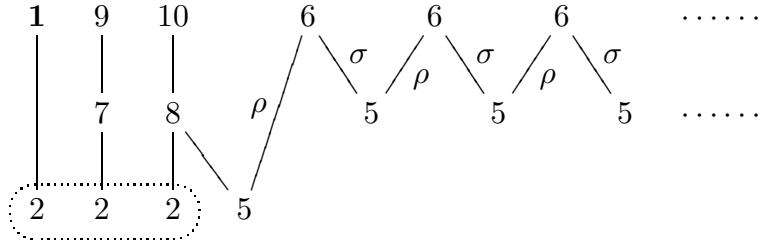
\ignore{
$$\xymatrixcolsep{0.75pc}\xymatrixrowsep{1.5pc}
\xymatrix{ {\bold 1}\edge[dd] &9\edge[d] &10\edge[d] &&6\edge[ddl]_\rho
\edge[dr]^\sigma &&6\edge[dl]^\rho \edge[dr]^\sigma &&6\edge[dl]^\rho
\edge[dr]^\sigma &&\cdots\cdots\\
 &7\edge[d] &8\edge[d] \edge[dr] &&&5 &&5 &&5 &\cdots\cdots\\ 2\threepool &2
&2 &5 }$$ }

\noindent This phantom is effective relative to the class of modules in
$\pinflamod$ obtained by chopping suitable `infinite tails' off the given
graph.

An example of an infinite dimensional $\pinflamod$-phan\-tom of $S_8$,
finally, is

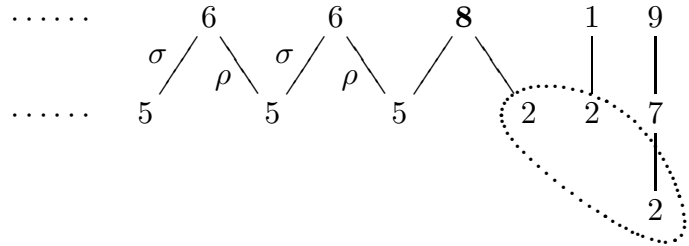
\ignore{
$$\xymatrixcolsep{0.75pc}\xymatrixrowsep{1.5pc}
\xymatrix{
\cdots\cdots &&6\edge[dl]_\sigma \edge[dr]_\rho &&6\edge[dl]_\sigma
\edge[dr]_\rho &&{\bold 8}\edge[dl] \edge[dr] &&1\edge[d] &9\edge[d]\\
\cdots\cdots &5 &&5 &&5 &&2
 \save[0,0]+(-3,0);[0,0]+(-3,0) **\crv{~*=<2.5pt>{.} [0,0]+(-3,4) 
&[0,1]+(4,4) &[1,2]+(4,0) &[1,2]+(3,-4) &[1,2]+(-3,-3) &[0,0]+(-3,-3)}\restore
&2 &7\edge[d]\\
 & && && && & &2 }$$ }
\enddefinition

This concludes our discussion of Examples A,B,C.  The following criterion
from \cite{\HaHZ} is hovering in the background of most of the infinite
dimensional phantoms displayed so far.

\definition{Criterion for failure of contravariant finiteness of
$\A$}

Let  $\la=\kgami$ be a split finite dimensional algebra, and, once more,
 suppose that $\A\subseteq\lamod$ is closed under finite direct sums. The
simple module $S=\la e_1/Je_1$ centered in the vertex $e_1$ fails to have an
$\A$-ap\-prox\-i\-ma\-tion in case the following holds:

The vertex $e_1$ can be supplemented to a sequence $e_1,\dots,e_m$ of distinct
vertices of $\Gamma$, together with sequences $p_1,\dots,p_m$, $q_1,\dots,q_m$
in $J$, where $p_i=p_ie_i$ and $q_i=q_ie_i$ are such that  conditions (1) and
(2) below are satisfied:

(1) For each $n\in\NN$, there exists a module $M_n\in\A$ having a graph that
contains a subgraph of the form

\ignore{
$$\xymatrixcolsep{0.4pc}\xymatrixrowsep{2.0pc}
\xymatrix{ e_1\edge[dr]_(0.6){p_1} &&e_2\edge[dl]_(0.4){q_2}
\edge[dr]_(0.6){p_2} &&e_3\edge[dl]_(0.4){q_3} \edge[dr]_(0.6){p_3}
&\ar@{.}[r] &&e_m\edge[dl]_(0.4){q_m} \edge[dr]_(0.6){p_m}
&&e_1\save+<0pc,1.5pc>\drop{n\hbox{\ repetitions}}\restore
\save[0,-9]+(0,-3);[0,9]+(3,4)**\frm{^)}\restore
\edge[dl]_(0.4){q_1}
\edge[dr]_(0.6){p_1} &\ar@{.}[r] &&e_1\edge[dl]_(0.4){q_1}
\edge[dr]_(0.6){p_1} &&e_2\edge[dl]_(0.4){q_2} \edge[dr]_(0.6){p_2}
&\ar@{.}[r] &&e_m\edge[dl]_(0.4){q_m} \edge[dr]_(0.6){p_m} &\\
 &\bullet &&\bullet &&\bullet\ar@{.}[r] &\bullet &&\bullet &&\bullet\ar@{.}[r]
&\bullet &&\bullet &&\bullet\ar@{.}[r] &\bullet &&\bullet }$$ }

(2) Given any object $A$ in $\A$, the top elements of $A$ of type $e_1$ are
not annihilated by $p_1$, and 
$$\bigl( p_ia=q_{i+1}b\ne0 \quad\implies\quad p_{i+1}b\ne0 \bigr)$$ for
$a,b\in A$ and $1\le i\le m$; here $p_{m+1}=p_1$ and $q_{m+1}=q_1$. \qed
\enddefinition

In fact, the hypotheses of the criterion yield an infinite dimensional factor
module of $\varinjlim M_n$ which is an $\A$-phan\-tom of $S$; it has a graph
containing the `infinite extension to the right' of the graph in condition
(1) as a subgraph.

\head 5. Phantoms over string algebras\endhead

In this section, we present a preview of ongoing joint work with S. O.
Smal\o\.

A special class of string algebras $\la_{m,n}$, for $m,n \ge 2$  -- those on
the quiver 

\ignore{
$$\xymatrix{ \bullet \ar@(ul,dl)_\alpha \ar@(ur,dr)^\beta }$$ }

\noindent  and subject to the relations $\alpha\beta =\beta\alpha =0$ and 
$\alpha^m= \beta^n =0$  --  was first singled out by Gelfand and Ponomarev in
the late 1960's, as being intimately related to the representation theory of
the Lorentz group \cite{\GePo}; in fact, classifying the finitely generated
indecomposable modules over the algebras $\la_{m,n}$ amounts to a
classification of the Harish-Chandra modules over the Lorentz group. Taking
this route, Gelfand and Ponomarev gave a hands-on structural description of
the finitely generated indecomposable objects in
$\la_{m,n}$-mod. In particular, their findings show that over an
algebraically closed base field the algebras
$\la_{m,n}$ are tame. Actually, this work was preceded by an investigation of
Szekeres \cite{\Sze} into the structure of the finitely generated modules
over certain factor rings of $\ZZ[x]$, which anticipates many of the ideas
re-encountered in later work on biserial algebras. In the mid-seventies,
Gabriel presented a categorical reinterpretation of the Gelfand-Ponomarev
approach (see
\cite{\Gab}), which in turn caused Ringel to recognize that these methods are 
applicable in a far wider context: In a first round of generalizations, he
used them to describe the finite dimensional indecomposable representations of
the dihedral 2-groups in characteristic 2
\cite{\Rin}; this work appeared in tandem with a paper of Bondarenko
containing roughly the same information \cite{\Bon}. Next, Donovan and
Freislich applied these methods to the `biserial' algebras introduced below 
\cite{\DoFr}. Clearly, the algebras $\la_{m,n}$ considered by Gelfand and
Ponomarev belong to the subclass of `special beserial' algebras; in fact,
they are even `string algebras'.

\definition{Definitions} (see \cite{\WaWa} and \cite{\BuRi}) (1) $\la$ is
called {\it biserial\/} if each indecomposable projective left or right
$\la$-module
$P$ has the following property: $\text{rad}\,P= U+V$, where $U,V$ are
uniserial (possibly trivial) with $U\cap V$ either zero or simple.

(2) $\la$ is {\it special biserial\/} provided  $\la$ is of the form
$\kgami$ such that 

$\bullet$ Given any vertex $e$ of $\Gamma$, there are at most two arrows
entering $e$ and at most two arrows leaving $e$, and

$\bullet$ Given any arrow $\alpha$ of $\Gamma$, there exists at most one
arrow $\beta$ such that 
$\beta\alpha$ does not belong to $I$, and at most one arrow $\gamma$ such that
$\alpha\gamma$ does not belong to $I$.

Moreover, $\la$ is a {\it string algebra\/} if, in addition, $\la$ is a
monomial relation algebra, meaning that $I$ can be generated by certain paths
in
$\Gamma$. 
\enddefinition

Clearly, special biserial algebras are biserial. All finite dimensional
biserial algebras over algebraically closed fields are known to be tame: the
special biserial case was completed by Wald and Waschb\"usch in \cite{\WaWa},
while the general biserial situation  was settled much later by
Crawley-Boevey \cite{\CBser} on the basis of an alternate description of
biserial algebras due to him and Vila-Freyer
\cite{\VFCB} and a remarkable result of Geiss \cite{\Geitame} (saying that
algebras with tame degenerations are always tame).

All the while, special biserial algebras have continued to provide challenges
which, in spite of the availability of a highly explicit classification of the
indecomposables, were far from resolvable at a glance. We mention only a few
such lines, together with a selection of references, not aiming at
completeness:

Let $\la$ be biserial.

$\bullet$ When does $\la$ have finite representation type? (See \cite{\SkWa}.)

$\bullet$ What does the Auslander-Reiten quiver of $\la$ look like? (See,
e.g., \cite{\BuRi}, \cite{\ErSk}, and \cite{\Geiser}.)

$\bullet$ Characterize the auto-equivalences of the category $\lamod$. (See
\cite{\Bleone} and \cite{\Bletwo}.)

$\bullet$ Describle the generic modules over string algebras and the maps
among them. (See
\cite{\Ringen}).

New sources of symmetric biserial algebras can be found in \cite{\Rog}.

Here we will provide an overview of a complete, constructive solution to the
problem as to which string algebras $\la$ have the property that $\pinflamod$
is contravariantly finite. Our proof of the answer -- not given here -- makes
full use of the description of the indecomposable objects in $\lamod$, as
reviewed graphically below.

\proclaim{Theorem 8} {\rm (Its evolution can be traced in \cite{\GePo},
\cite{\Rin},
\cite{\DoFr}, \cite{\WaWa})} Let $\la=\kgami$ be a special biserial algebra
over an algebraically closed field $K$. Then each indecomposable  object in
$\lamod$ is either a band module or a string module. Here the string modules
are those with graphs of the form

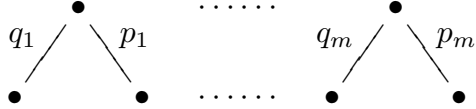
\ignore{
$$\xymatrixcolsep{0.75pc}\xymatrixrowsep{1.5pc}
\xymatrix{
 &\bullet \edge[dl]_{q_1} \edge[dr]^{p_1} &&\cdots\cdots &&\bullet
\edge[dl]_{q_m} \edge[dr]^{p_m}\\
\bullet &&\bullet &\cdots\cdots &\bullet &&\bullet }$$ }

\noindent where the $p_i,q_i$ are paths in $K\Gamma\setminus I$, with $q_1$
and $p_m$ possibly trivial, such that $\operatorname{first arrow}(q_i) \ne
\operatorname{first arrow}(p_i)$ for $1 \le i \le m$, and
$\operatorname{last arrow}(p_i)\ne \operatorname{last arrow}(q_{i+1})$ for
$1\le i\le m-1$. (Note that some of these conditions may be void, namely in
case the relevant paths are trivial.)

The band modules are characterized by their graphs, paired with irreducible
vector space automorphisms, as follows: The pertinent graphs are of the form

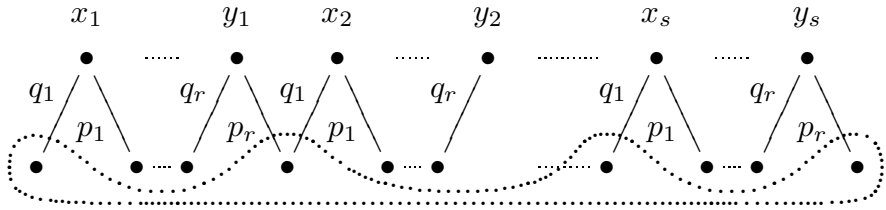
\ignore{
$$\xymatrixcolsep{0.4pc}\xymatrixrowsep{2.0pc}
\xymatrix{
 &\bullet\toplabel{x_1} \edge[dl]_(0.4){q_1} \edge[dr]_(0.6){p_1} &\ar@{.}[r]
&&\bullet\toplabel{y_1} \edge[dl]_(0.4){q_r} \edge[dr]_(0.6){p_r}
&&\bullet\toplabel{x_2} \edge[dl]_(0.4){q_1} \edge[dr]_(0.6){p_1} &\ar@{.}[r]
&&\bullet\toplabel{y_2} \edge[dl]_(0.4){q_r} &\ar@{.}[rr]
&&&\bullet\toplabel{x_s} \edge[dl]_(0.4){q_1} \edge[dr]_(0.6){p_1} &\ar@{.}[r]
&&\bullet\toplabel{y_s} \edge[dl]_(0.4){q_r} \edge[dr]_(0.6){p_r}\\
\bullet
 \save[0,0]+(-3,0);[0,0]+(-3,0) **\crv{~*=<2.5pt>{.} [0,0]+(-3,4) 
&[0,0]+(3,4) &[0,2]+(-3,-3) &[0,3]+(3,-3) &[0,5]+(-3,4) &[0,5]+(3,4)
&[0,7]+(-3,-3) &[0,10]+(3,-3) &[0,12]+(-3,4) &[0,12]+(3,4) &[0,14]+(-3,-3)
&[0,15]+(3,-3) &[0,17]+(-3,4) &[0,17]+(3,4) &[0,17]+(3,-4) &[0,17]+(-3,-4)
&[0,0]+(3,-4) &[0,0]+(-3,-4)}\restore
 &&\bullet\ar@{.}[r] &\bullet &&\bullet &&\bullet\ar@{.}[r] &\bullet
&&\ar@{.}[rr] &&\bullet &&\bullet\ar@{.}[r] &\bullet &&\bullet }$$  }

\noindent where $p_i,q_i$ are are nontrivial paths in $K\Gamma\setminus I$
with $\operatorname{first arrow}(q_i)\ne \operatorname{first arrow}(p_i)$ for
all $i$ and $\operatorname{last arrow}(p_i)\ne \operatorname{last
arrow}(q_{i+1})$ for $i<r$, and also $\operatorname{last arrow}(p_r)\ne
\operatorname{last arrow}(q_1)$. The nature of the dotted pool is specified by
the dependence relation
$$p_ry= \sum_{i=1}^s k_iq_1x_i,$$ where 
$\left[ \smallmatrix 0 &\cdots &0 &k_1\\ 1 &\cdots &0 &k_2\\
\hphantom{\cdots} &\ddots &\hphantom{\cdots} &\vdots\\  0 &\cdots &1
&k_s\endsmallmatrix \right]$
 is the Frobenius companion matrix of an irreducible automorphism of $K^s$.

Moreover, all modules having one of the above descriptions are indecomposable.
\qed\endproclaim

This classification can be completed with a suitable uniqueness statement
which does not impinge on our present results.  For our main theorem,  string
modules will be of particular relevance. A  {\it generalized string module\/}
will be a module
$X\in\Lamod$ which arises as a direct limit of a countable directed system of
string modules, each embedded into its successor, i.e., any module $X$ having
a graph of one of the following forms

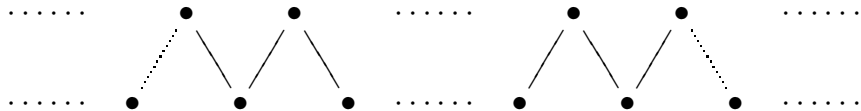
\ignore{
$$\xymatrixcolsep{0.5pc}\xymatrixrowsep{1.5pc}
\xymatrix{
\cdots\cdots &&\bullet\ar@{.}[dl] \edge[dr] &&\bullet\edge[dl] \edge[dr]
&&\cdots\cdots &&\bullet\edge[dl] \edge[dr] &&\bullet\edge[dl] \ar@{.}[dr]
&&\cdots\cdots\\
\cdots\cdots &\bullet &&\bullet &&\bullet &\cdots\cdots &\bullet &&\bullet
&&\bullet &\cdots\cdots }$$ }

\noindent such that each finite segment is the graph of a string module.
(Here the dotted edges may but need not appear, depending on whether the
graph is one-/two-sided infinite or finite.)  These modules were also
considered by Ringel in \cite{\Rintwo} as ``modules associated with
$\NN$-words of $\ZZ$-words'' over the alphabet
$\Gamma_0\cup
\Gamma_0^{-1}$, where $\Gamma_0$ is the vertex set of the quiver $\Gamma$ of
$\la$.  It should be self-explanatory what we mean by a {\it left and right
periodic\/} generalized string module; we merely emphasize that we allow 
the ``left'' and ``right'' periods to differ and that we regard termination
as a period.

We are now in a position to state the main new result of this section.

\proclaim{Theorem 9} {\rm \cite{\HZSmbsrl}} Let $\la=\kgami$ be a finite
dimensional string algebra.  For any simple $S\in\lamod$, there exists a
generalized string module $H=H(S)$ which is uniquely determined by $\Gamma$
and $I$, together with a canonical homomorphism $f : H\rightarrow S$, having
the following properties:
\smallskip

\noindent{\rm (I)} The following statements are equivalent:

{\rm (i)} $\dim_K H<\infty$;

{\rm (ii)} $S$ has a classical \pinfapprox;

{\rm (iii)} $f : H\rightarrow S$ is the minimal \pinfapprox\ of $S$.
\smallskip

\noindent{\rm (II)} $H$ always belongs to $\pinf(\Lamod)$ and is a
$\pinflamod$-phan\-tom of
$S$.  In fact, the map $f : H\rightarrow S$ makes $H$ an effective
$\pinflamod$-phan\-tom with respect to the class $\sinflamod$ of all string
modules of finite projective dimension.
\smallskip

\noindent{\rm (III)} $H$ is left and right periodic and can be constructed
from
$\Gamma$ and $I$ in fewer than $3|\Gamma_0|$ steps.
\qed\endproclaim

Given a simple module $S=\la e/Je$ over a string algebra $\la$, we will call
the module
$H=H(S)$ of the theorem the {\it characteristic\/}
$\pinflamod$-phan\-tom of $S$.  After illustrating the theory with a first
example, we will give an inductive description of the finite segments of a
graph

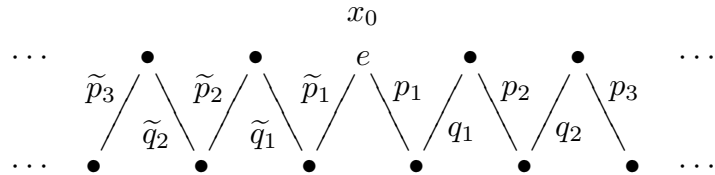
\ignore{
$$\xymatrixcolsep{0.5pc}\xymatrixrowsep{2.0pc}
\xymatrix{
\cdots &&\bullet \edge[dl]_(0.4){\ptil_3} \edge[dr]_(0.6){\qtil_2} &&\bullet
\edge[dl]_(0.4){\ptil_2} \edge[dr]_(0.6){\qtil_1} &&e\toplabel{x_0}
\edge[dl]_(0.4){\ptil_1}
\edge[dr]^(0.4){p_1} &&\bullet \edge[dl]^(0.6){q_1}
\edge[dr]^(0.4){p_2} &&\bullet \edge[dl]^(0.6){q_2}
\edge[dr]^(0.4){p_3} &&\cdots\\
\cdots &\bullet &&\bullet &&\bullet &&\bullet &&\bullet &&\bullet &\cdots }$$
}

\noindent of the characteristic $\pinf(\lamod)$-phantom of $S$. Sending the
top element marked $x_0$ to
$e+Je$ in $S$ and sending the other top elements to zero will then yield a
map $f : H\rightarrow S$ as stipulated in the theorem.

\definition{Example F} Let $\la=\kgami$ be the string algebra with the
following indecomposable projective left $\la$-modules

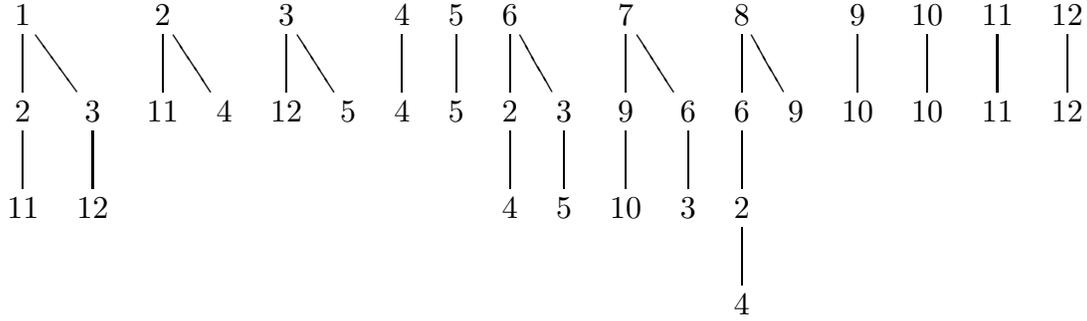
\ignore{
$$\xymatrixcolsep{0.5pc}\xymatrixrowsep{1.5pc}
\xymatrix{ 1 \edge[d] \edge[dr] &&2 \edge[d] \edge[dr] &&3 \edge[d] \edge[dr]
&&4
\edge[d] &5 \edge[d] &6 \edge[d] \edge[dr] &&7 \edge[d] \edge[dr] &&8
\edge[d] \edge[dr] &&9 \edge[d] &10 \edge[d] &11 \edge[d] &12 \edge[d]\\ 2
\edge[d] &3 \edge[d] &11 &4 &12 &5 &4 &5 &2 \edge[d] &3 \edge[d] &9
\edge[d] &6 \edge[d] &6 \edge[d] &9 &10 &10 &11 &12\\  11 &12 && && & & &4 &5
&10 &3 &2 \edge[d]\\ & && && & & && && &4 }$$ }

\noindent Then the characteristic $\pinflamod$-phan\-tom $H_1$ of $S_1$ has
graph

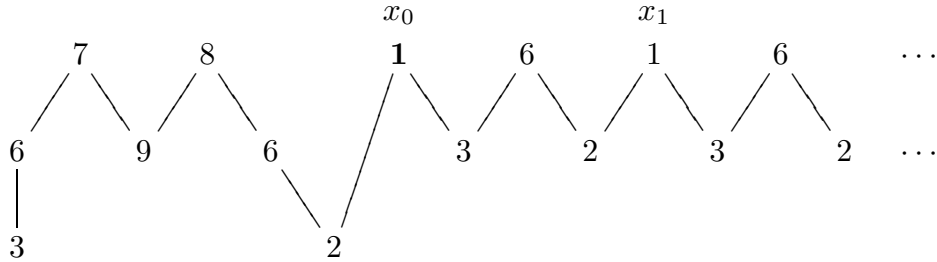
\ignore{
$$\xymatrixcolsep{0.75pc}\xymatrixrowsep{1.5pc}
\xymatrix{
 &7 \edge[dl] \edge[dr] &&8 \edge[dl] \edge[dr] &&&{\bold 1}\toplabel{x_0}
\edge[ddl]
\edge[dr] &&6 \edge[dl] \edge[dr] &&1\toplabel{x_1} \edge[dl]
\edge[dr] &&6
\edge[dl] \edge[dr]  &&\cdots\\ 6 \edge[d] &&9 &&6 \edge[dr] &&&3 &&2 &&3 &&2
&\cdots\\ 3 &&&&&2 }$$ }

\noindent and the homomorphism $f : H_1\rightarrow S_1$ which sends $x_0$ to
a nonzero element of
$S_1$ and the $x_i$ for $i\ge1$ to zero has the properties described in
Theorem 9. In particular,
$S_1$ does not have a
\pinfapprox.

The characteristic $\pinflamod$-phan\-tom of $S_7$, on the other hand, has
graph

\ignore{
$$\xymatrixrowsep{1.25pc}
\xymatrix{ 7\edge[d]\\ 6\edge[d]\\ 3}$$ }

\noindent and thus coincides with the minimal \pinfapprox\ of $S_7$. These
graphs are obtained via the algorithm which we describe next.
\qed\enddefinition

Recall that we refer to a module $X$ as a top-embeddable submodule of $Y$ if
there exists a monomorphism $f : X\rightarrow Y$ which induces a monomorphism
$X/JX\rightarrow Y/JY$. Dually, we call $X$ a {\it socle-faithful\/} factor
module of
$Y$ if there exists an epimorphism $f : Y\rightarrow X$ which induces an
epimorphism $\soc Y\rightarrow \soc X$.
\medskip

\noindent{\bf Description of the characteristic phantom of a simple module}
$S\in\lamod$, where $\la$ is a finite dimensional string algebra.

Let $S=\la e/Je$. The following are the steps of an algorithmic procedure for
constructing $H=H(S)$, but here we will not discuss the algorithmic nature,
nor prove that the quantities stipulated in the process exist.

\underbar{Step 1}. Let $p_1$ and $\ptil_1$ be paths starting in $e$ which have
{\it minimal\/} lengths $\ge0$ such that

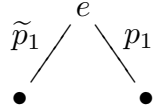
\ignore{
$$\xymatrixcolsep{0.75pc}\xymatrixrowsep{1.5pc}
\xymatrix{
 &e \edge[dl]_{\ptil_1} \edge[dr]^{p_1}\\
\bullet &&\bullet }$$ }

\noindent is the graph of a string module which can be top-embedded into an
object in $\sinflamod$. In particular, we have $p_1\ne \ptil_1$ unless both of
these paths are trivial, and $\text{startarrow}(p_1)\ne
\text{startarrow}(\ptil_1)$ if both are nontrivial.

If both $p_1$ and $\ptil_1$ are trivial, we set $H=S$. Otherwise, we proceed
to

\underbar{Step 2}. Let $q_1$ and $\qtil_1$ be paths ending in
$\text{end}(p_1)$ and $\text{end}(\ptil_1)$, respectively, and having {\it
maximal\/} lengths $\ge0$ with the property that

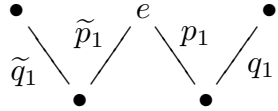
\ignore{
$$\xymatrixcolsep{0.75pc}\xymatrixrowsep{1.5pc}
\xymatrix{
\bullet \edge[dr]_{\qtil_1} &&e \edge[dl]_(0.45){\ptil_1}
\edge[dr]^(0.45){p_1} &&\bullet \edge[dl]^{q_1}\\
 &\bullet &&\bullet }$$ }

\noindent is the graph of a string module which arises as a socle-faithful
factor module of an object in $\sinflamod$. In case $p_1$ is trivial, set
$q_1=e$, and deal symmetrically with $\qtil_1$.

If both $q_1$ and $\qtil_1$ are trivial, i.e., if the graph of Step 2
coincides with that of Step 1, we let $H$ be the string module having this
graph. Otherwise, we proceed to

\underbar{Step 3}. Let $p_2$ and $\ptil_2$ be paths starting in
$\text{start}(q_1)$ and $\text{start}(\qtil_1)$, respectively, which have {\it
minimal\/} lengths $\ge0$ with the property that 

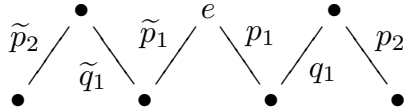
\ignore{
$$\xymatrixcolsep{0.75pc}\xymatrixrowsep{1.5pc}
\xymatrix{
 &\bullet \edge[dl]_{\ptil_2} \edge[dr]_(0.55){\qtil_1} &&e
\edge[dl]_(0.45){\ptil_1} \edge[dr]^(0.45){p_1} &&\bullet
\edge[dl]^(0.55){q_1} \edge[dr]^{p_2}\\
\bullet &&\bullet &&\bullet &&\bullet }$$ }

\noindent is the graph of a string module which can be top-embedded into some
object in $\sinflamod$, \dots etc.

After fewer than $3|\Gamma_0|$ steps, this procedure either terminates or has
become periodic on both sides. It is easy to recognize when one has hit a left
or right period: Indeed, if $\text{startarrow}(p_i)= \text{startarrow}(p_j)$
for some $i<j$, then $p_{i+r}=p_{j+r}$ and $q_{i+r}=q_{j+r}$ for all $r\ge0$,
the symmetric criterion holding for the other side. \qed
\medskip

As one gleans from this inductive description of the phantoms $H=H(S)$ of the
simple modules $S$, string algebras merit their name also from a homological
viewpoint. Indeed, the phantoms $H(S)$ depend only on the string modules of
finite projective dimension, and the $K$-dimensions of the $H(S)$ in turn
determine whether or not $\pinflamod$ is contravariantly finite. One may
wonder whether our emphasis on string modules is just dictated by
convenience and whether band modules play a similar homological role.  They
do not.  In fact, there exist string algebras of positive little finitistic
dimension which have no nontrivial (finitely generated) band modules of
finite projective dimension. Another asset of string algebras that arises as a
byproduct of our main theorem we record somewhat more formally.

\proclaim{Corollary 10} Suppose that $S$ is a simple module over a finite
dimensional string algebra $\la$. If $S$ has a \pinfapprox, then the minimal
\pinfapprox\ of $S$ is a string module and, in particular, is indecomposable.
\qed\endproclaim

The corollary does not carry over to arbitrary special biserial algebras, as
the next example demonstrates.

\definition{Example G} Let $\la=\kgami$ be the special biserial algebra with
indecomposable projectives

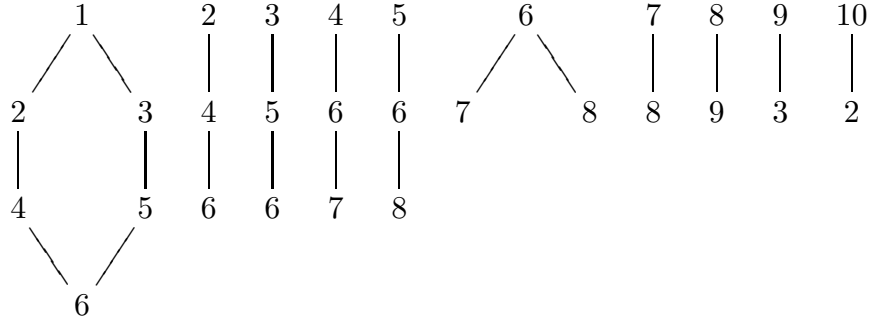
\ignore{
$$\xymatrixcolsep{0.75pc}\xymatrixrowsep{1.5pc}
\xymatrix{
 &1 \edge[dl] \edge[dr] &&2 \edge[d] &3 \edge[d] &4 \edge[d] &5 \edge[d] &&6
\edge[dl] \edge[dr] &&7 \edge[d] &8 \edge[d] &9\edge[d] &10\edge[d]\\ 2
\edge[d] &&3 \edge[d] &4 \edge[d] &5 \edge[d] &6 \edge[d] &6 \edge[d] &7 &&8
&8 &9 &3 &2\\ 4 \edge[dr] &&5 \edge[dl] &6 &6 &7 &8\\
 &6 }$$ }

\noindent Then the minimal \pinfapprox\ of $S_1$ is as follows:

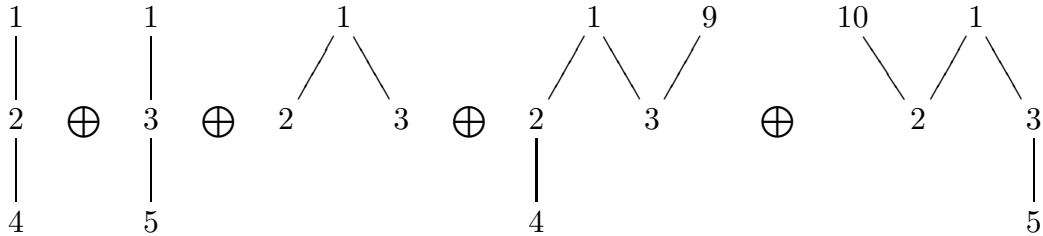
\ignore{
$$\xymatrixcolsep{0.6pc}\xymatrixrowsep{1.5pc}
\xymatrix{ 1 \edge[d] &&1 \edge[d] &&&1 \edge[dl] \edge[dr] &&&&1\edge[dl]
\edge[dr] &&9\edge[dl] &&10\edge[dr] &&1\edge[dl] \edge[dr]\\ 2 \edge[d]
&\bigoplus &3 \edge[d] &\bigoplus &2 &&3 &\bigoplus &2\edge[d] &&3
&&\bigoplus &&2 &&3\edge[d]\\ 4 &&5 &&&&&&4 &&&& &&&&5 }$$ }

\noindent In particular, we observe that the factor modules of $\la e_1$
contained in the two rightmost summands are not minimal with respect to
top-embeddability into objects of $\pinflamod$. Thus the alternation ``choose
a minimal factor module embeddable into a module in $\pinflamod$, then choose
a maximal essential extension arising as a factor module of a module in
$\pinflamod$'' which leads to the minimal approximations of the simples over
string algebras in case of existence, cannot be expected to achieve this goal
in the more general situation.
\qed\enddefinition

To give another illustration of our algorithm, we conclude with an example of
a string algebra $\la$ and a simple left
$\la$-module $S$, the characteristic phantom of which is twosided infinite
with  left/right periods reached at different steps of the algorithm.

\definition{Example H} Let $\la=\kgami$ be the string algebra with the
following indecomposable projective left modules:

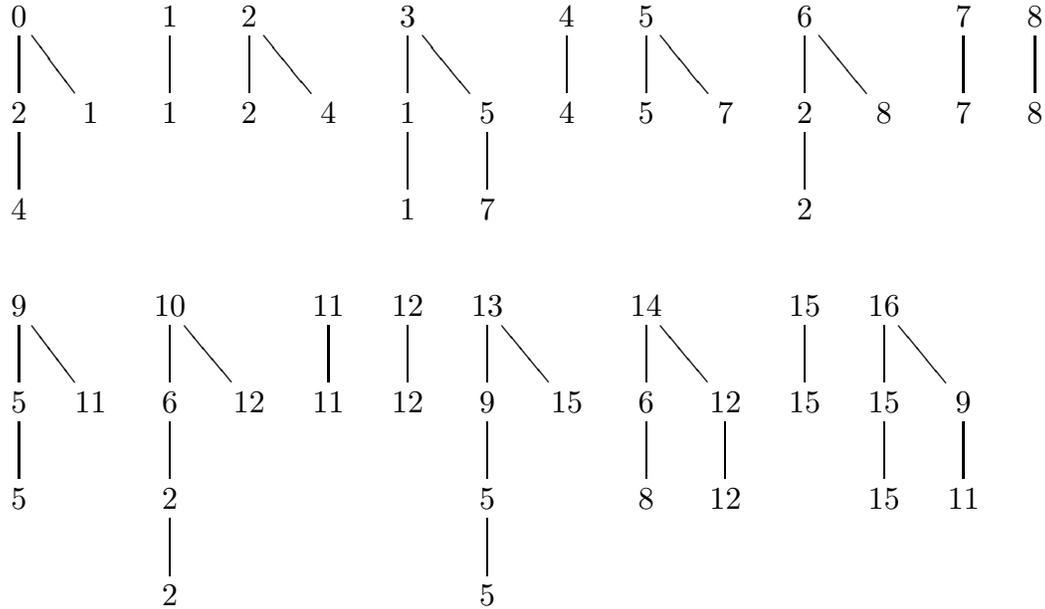
\ignore{
$$\xymatrixcolsep{0.75pc}\xymatrixrowsep{1.5pc}
\xymatrix{ 0\edge[d] \edge[dr] &&1\edge[d] &2\edge[d] \edge[dr] &&3\edge[d]
\edge[dr] &&4\edge[d] &5\edge[d] \edge[dr] &&6\edge[d] \edge[dr] &&7\edge[d]
&8\edge[d]\\ 2\edge[d] &1 &1 &2 &4 &1\edge[d] &5\edge[d] &4 &5 &7 &2\edge[d]
&8 &7 &8\\ 4 & &&&&1 &7 &&&&2\\ 9\edge[d] \edge[dr] &&10\edge[d] \edge[dr]
&&11\edge[d] &12\edge[d] &13\edge[d] \edge[dr] &&14\edge[d] \edge[dr]
&&15\edge[d] &16\edge[d]
\edge[dr]\\ 5\edge[d] &11 &6\edge[d] &12 &11 &12 &9\edge[d] &15 &6\edge[d]
&12\edge[d] &15 &15\edge[d] &9\edge[d]\\ 5 &&2\edge[d] &&&&5\edge[d] &&8 &12
&&15 &11\\
 &&2 &&&&5 }$$ }

\noindent Then the characteristic $\pinflamod$-phan\-tom of $S= \la e_0/Je_0$
has a graph as follows:

\ignore{
$$\xymatrixcolsep{0.6pc}\xymatrixrowsep{1.5pc}
\xymatrix{
\cdots &&14\edge[dl]\edge[d] &10\edge[dl]\edge[d] &14\edge[dl]\edge[d]
&10\edge[dl]\edge[d] &{\bold 0}\edge[dd]\edge[dr] &&3\edge[dl]\edge[dd]
&13\edge[d]\edge[dr] &16\edge[d]\edge[dr] &13\edge[d]\edge[dr]
&16\edge[d]\edge[dr] &&\cdots\\
\cdots &6 &12 &6 &12 &6\edge[dr] &&1 &&9\edge[dl] &15 &9 &15 &9 &\cdots\\
 &&&&&&2 &&5 }$$ }
\enddefinition

\enddocument